\documentclass[12pt]{article}

\usepackage{amssymb,latexsym,amsmath,setspace}
\usepackage{amsthm,amsfonts}
\usepackage{graphicx}
\usepackage [left =1.0in, top =1.0in, bottom =1.0in,
    right=1.0in]{geometry}

\newcommand{\abs}[1]{\left\vert#1\right\vert}
\newcommand{\Lip}[1]{\text {Lip}\left(#1\right)}

\newcommand{\C}{cnst}

\theoremstyle{plain}
\newtheorem{thm}{Theorem}[section]
\newtheorem{cor}[thm]{Corollary}
\newtheorem{lem}[thm]{Lemma}

\theoremstyle{definition}

\theoremstyle{remark}
\newtheorem{rem}{Remark}[section]

\newcommand{\comment}[1]{}
\begin{document}

\title {\textbf{A simple piston problem in one dimension}
\\(accepted by Nonlinearity\footnote {
Nonlinearity \copyright 2006 IOP Publishing Ltd.
\texttt{http://www.iop.org/EJ/journal/Non}}~)}

\author {Paul Wright\\Department of Mathematics, New York University\\
        251 Mercer St.\\
        New York, NY 10012\\\texttt{paulrite@cims.nyu.edu}}

\date {August 2006}
\maketitle

\begin {abstract}

We study a heavy piston that separates finitely many ideal gas
particles moving inside a one-dimensional gas chamber.  Using
averaging techniques, we prove precise rates of convergence of the
actual motions of the piston to its averaged behavior. The
convergence is uniform over all initial conditions in a compact set.
The results extend earlier work by Sinai and Neishtadt, who
determined that the averaged behavior is periodic oscillation. In
addition, we investigate the piston system when the particle
interactions have been smoothed.  The convergence to the averaged
behavior again takes place uniformly, both over initial conditions
and over the amount of smoothing.

\end {abstract}




\section{Introduction} \label {sct:intro}

Consider the following simple model of a piston separating two gas
chambers:  A piston of mass $M\gg 1$ divides a cylindrical chamber
in $\mathbb{R}^3$ into two halves. The piston is parallel to the two
ends of the chamber and can only move in the normal direction. On
either side of the piston there are a finite number of gas particles
of unit mass. All of the gas particles are point particles that
interact with the walls of the chamber and with the piston via
elastic collisions. The interactions of the gas particles with the
piston and the ends of the chamber are completely specified by their
motions along the normal axis of the chamber. Thus, this system
projects onto a system inside the unit interval consisting of a
massive point particle, the piston, which interacts with the gas
particles on either side of it. These gas particles make elastic
collisions with the walls at the ends of the chamber and with the
piston, but they do not interact with each other:  They form an
ideal gas.

This simple model is useful for investigating the adiabatic piston
problem~\cite {Ca63}.  This well-known problem from physics concerns
an adiabatic piston, i.e.~one with no internal degrees of freedom,
that separates two gas chambers. Initially, the piston is fixed in
place, and the gas in each chamber is in a separate thermal
equilibrium. At some time, the piston is no longer externally
constrained and is free to move. One hopes for an ergodic theorem to
show that eventually the system will come to a full thermal
equilibrium, where each gas has the same pressure and temperature.
If one existed, then the final pressure and temperature could be
predicted using thermodynamics.  However, whether the system will
evolve to thermal equilibrium and the interim behavior of the piston
are mechanical problems that are not yet explicitly
resolved~\cite{Gru99}.  We do not give a full investigation of these
issues here, but we wish to emphasize that they cannot be addressed
by thermodynamics alone. The interested reader should see some
papers~\cite {Che05,GPL03} and the references therein for details
about recent progress in this area.

Gruber \emph{et al.}~\cite{GF99,GPL02,GPL03} have extensively
studied the simple model above using the Boltzmann equation, the
Liouville equation, and numerical simulations. Making various
assumptions, they observed that the system evolves in at least two
stages. First, the system relaxes deterministically and
adiabatically toward mechanical equilibrium, where the pressures on
either side of the piston are equal.  In the second, much longer,
stage, the piston drifts stochastically in the direction of the
hotter gas, and the temperatures of the gases equilibrate.  Because
of the piston's stochastic fluctuations, heat is allowed to flow
between the gases, and in this sense the piston is no longer
adiabatic. It should be emphasized that the behavior of the piston
depends strongly on the ratios $n_i/M$, where $n_1$ and $n_2$ are
the number of gas particles on either side of the piston. In the
thermodynamic limit $M,n_1,n_2\rightarrow\infty$ while $n_1/M$ and
$n_2/M$ are held fixed, Gruber \emph{et al.}~\cite{GPL02} concluded
that in the first stage above, the piston performs damped
oscillatory motion, where the damping is strong if $n_1/M,n_2/M>1$,
and weak if $n_1/M,n_2/M<1$. Using kinetic theory, Crosignani
\emph{et al.}~\cite{CD96} had already derived similar equations
describing this damped oscillatory motion for the adiabatic piston.

Sinai~\cite {Sin99} also investigated the simple model, but with
averaging techniques that examine the limit where
$M\rightarrow\infty   $ while the total energy of the system is
bounded and $n_1$ and $n_2$ are fixed. He determined that the
averaged behavior of the piston is periodic oscillation, with the
piston moving inside an effective potential well whose shape depends
on the initial position of the piston and the initial energies of
the gases. Neishtadt and Sinai~\cite{NS04} pointed out that a
classical averaging theorem due to Anosov, proved for smooth
systems, can be extended to this case. Their insight allows us to
conclude that if we examine the actual motions of the piston with
respect to the slow time $\tau=t/M^ {1/2}$ as $M\rightarrow\infty $,
then in probability (with respect to Riemannian volume) most initial
conditions give rise to orbits whose actual motion is accurately
described by the averaged behavior for $\tau\in [0, 1]$.

This paper proves that the actual motions do not deviate by more
than $\mathcal{O} (M^ {-1/2}) $ from the averaged behavior for
$\tau\in [0,1] $, i.e.~for $t\in [0,M^ {1/2} ]$. Furthermore, the
size of the deviations is bounded, independent of the initial
conditions.

We also investigate the behavior of the system when the interactions
of the gas particles with the walls and the piston have been
smoothed, so that Anosov's theorem applies directly.  Let
$\delta\geq 0$ be a parameter of smoothing, so that $\delta=0 $
corresponds to the hard core setting above.  Then the averaged
behavior of the piston is still a periodic oscillation, which
depends smoothly on $\delta$.  We show that the deviations of the
actual motions of the piston from the averaged behavior are again
not more than $\mathcal{O} (M^ {-1/2}) $ on the time scale $M^ {1/2}
$. The size of the deviations is bounded uniformly, both over
initial conditions and over the amount of smoothing.

\section{Statement of results}

\label {sct:results}

\subsection{The hard core piston problem}
\label {sct:hard_core_results}

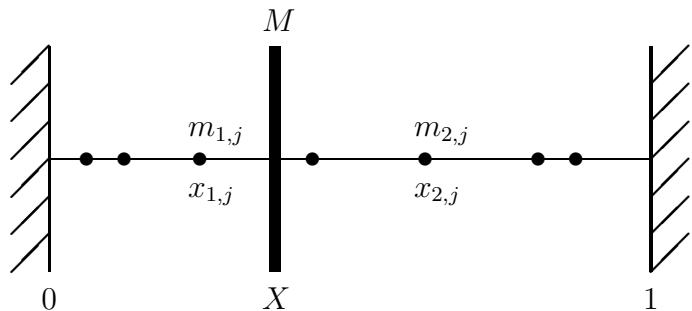
\begin{figure}
    \begin {center}
        \setlength{\unitlength}{1 cm}
    \begin{picture}(10,5)
        \thinlines
        \put(1,2.5){\line(1,0){8}}
        \thicklines
        \put(1,1){\line(0,1){3}}
        \put(9,1){\line(0,1){3}}
        \multiput(0.5,1)(0,0.5){6}{\line(1,1){0.5}}
        \multiput(9,1)(0,0.5){6}{\line(1,1){0.5}}
        \put(0.9,0.5){$0$}
        \put(8.9,0.5){$1$}
        \linethickness {0.15cm}
        \put(4.0,1){\line(0,1){3}}
        \put(3.83,0.5){$X$}
        \put(3.83,4.2){$M$}
        \put(1.5,2.5){\circle*{.18}}
        \put(2,2.5){\circle*{.18}}
        \put(3,2.5){\circle*{.18}}
        \put(4.5,2.5){\circle*{.18}}
        \put(6,2.5){\circle*{.18}}
        \put(7.5,2.5){\circle*{.18}}
        \put(8,2.5){\circle*{.18}}
        \put(2.85,2){$x_{1,j}$}
        \put(2.85,2.8){$m_{1,j}$}
        \put(5.85,2){$x_{2,j}$}
        \put(5.85,2.8){$m_{2,j}$}
    \end{picture}
    \end {center}
    \caption{The piston system with $n_1=3$ and $n_2=4$.  Note that the
        gas particles do not interact with each other, but only with the
        piston and the walls.}\label{fig:piston}
\end{figure}

Consider the system of $n_1+n_2+1$ point particles moving inside the
unit interval indicated in Figure \ref{fig:piston}.  One
distinguished particle, the piston, has position $X$ and mass $M$.
To the left of the piston there are $n_1>0$ particles with positions
$x_{1,j}$ and masses $m_{1,j} $, $1\leq j\leq n_1$, and to the right
there are $n_2>0$ particles with positions $x_{2,j}$ and masses
$m_{2,j} $, $1\leq j\leq n_2$.  These gas particles do not interact
with each other, but they interact with the piston and with walls
located at the end points of the unit interval via elastic
collisions.  We denote the velocities by $dX/dt=V$ and
$dx_{i,j}/dt=v_{i,j}$. There is a method for transforming this
system into a billiard system consisting of a point particle moving
inside an $(n_1+n_2+1)$-dimensional polytope~\cite{CM06}, but we
will not use this in what follows.

We are interested in the dynamics of this system when the numbers
and masses of the gas particles are fixed, the total energy is
bounded, and the mass of the piston tends to infinity. When
$M=\infty $, the piston remains at rest, and each gas particle
performs periodic motion.  More interesting are the motions of the
system when $M$ is very large but finite. Because the total energy
of the system is bounded, $MV^2/2\leq \C$, and so $V=\mathcal{O} (M^
{-1/2})$. Set $\varepsilon = M^ {-1/2} $, and let $W=V/\varepsilon
$, so that $dX/dt=\varepsilon W$ with $W =\mathcal{O} (1) $.

When $\varepsilon =0 $, the system has $n_1+n_2+2$ independent first
integrals (conserved quantities), which we take to be $X,\: W$, and
$s_{i,j} =\abs{v_{i,j}}$, the speeds of the gas particles. We refer
to these variables as the slow variables because they should change
slowly with time when $\varepsilon$ is small, and we denote them by
$h=(X,W, s_{1,j},s_{2,j})\in\mathbb{R}^ {n_1+n_2+2}$. Let
$h_\varepsilon(t,z) =h_\varepsilon(t)$ denote the dynamics of these
variables in time for a fixed value of $\varepsilon$, where $z$
represents the dependence on the initial condition in phase space.
We usually suppress the initial condition in our notation. Think of
$h_\varepsilon(\cdot) $ as a random variable which, given an initial
condition in the $2(n_1+n_2+1)$-dimensional phase space, produces a
piecewise continuous path in $\mathbb{R}^ {n_1+n_2+2}$. These paths
are the projection of the actual motions in our phase space onto a
lower dimensional space. The goal of averaging is to find a vector
field on $\mathbb{R}^ {n_1+n_2+2}$ whose orbits approximate
$h_\varepsilon(t) $.

\subsubsection {The averaged equation}
\label{sct:1dps_avg}

Sinai~\cite {Sin99} derived
\begin {equation} \label {eq:1davg}
    \frac{d}{d\tau}
    \begin {bmatrix}
    X\\
    W\\
    s_{1,j}\\
    s_{2,j}\\
    \end {bmatrix}
    =\bar H(h) :=
    \begin {bmatrix}
    W\\
    \frac{\sum_{j=1} ^ {n_1}m_{1,j} s_{1,j}^2}{X}-
    \frac{\sum_{j=1} ^ {n_2}m_{2,j} s_{2,j}^2}{1-X}\\
    -\frac{s_{1,j} W}{X}\\
    +\frac{s_{2,j} W}{1-X}\\
    \end {bmatrix}
\end {equation}
as the averaged equation (with respect to the slow time
$\tau=\varepsilon t$) for the slow variables.  We provide a
heuristic derivation in Section \ref{sct:1dps_hc_heuristic}.  Sinai
solved this equation as follows: From
$d\ln(s_{1,j})/d\tau=-d\ln(X)/d\tau$,
$s_{1,j}(\tau)=s_{1,j}(0)X(0)/X(\tau)$.  Similarly,
$s_{2,j}(\tau)=s_{2,j}(0) (1 - X(0))/(1 - X(\tau))$.  Hence
\[
    \frac {d^2X} {d\tau^2}=
    \frac{\sum_{j=1} ^ {n_1}m_{1,j} s_{1,j}(0)^2 X(0) ^2}{X^3}
    -\frac{\sum_{j=1} ^ {n_2}
    m_{2,j} s_{2,j}(0)^2 (1-X(0)) ^2}{(1-X) ^3},
\]
and so $(X, W) $ behave as if they were the coordinates of a
Hamiltonian system describing a particle undergoing periodic motion
inside a potential well. If we let $E_i=\sum_{j=1} ^ {n_i}m_{i,j}
s_{i,j}^2/2$ be the kinetic energy of the gas particles on one side
of the piston, the effective Hamiltonian may be expressed as
\begin {equation}
\label {eq:1dpot}
    \frac {1} {2}W^2+
    \frac{E_1(0)X(0)^2}{X^2}+\frac{E_2(0)(1-X(0))^2}{(1-X)^2}.
\end {equation}
Hence, the solutions to the averaged equation are periodic for all
initial conditions under consideration.

\subsubsection{Main result in the hard core setting} \label {sct:1dres}

The solutions of the averaged equation approximate the motions of
the slow variables, $h_\varepsilon(t) $, on a time scale
$\mathcal{O} (1/\varepsilon) $ as $\varepsilon\rightarrow 0$.
Precisely, let $\bar{h} (\tau,z)=\bar{h} (\tau) $ be the solution of
\[
\frac {d\bar{h}}{d\tau} =\bar {H} (\bar {h}),\qquad \bar {h} (0)
=h_\varepsilon(0).
\]
Again, think of $\bar h(\cdot) $ as being a random variable that
takes an initial condition in our phase space and produces a path in
$\mathbb R^{n_1+n_2+2}$.

Next, fix a compact set $\mathcal{V}\subset \mathbb R^{n_1+n_2+2}$
such that $h\in \mathcal{V} \Rightarrow X\subset\subset
(0,1),W\subset\subset \mathbb R$, and $s_{i,j} \subset\subset
(0,\infty)$ for each $i$ and $j$.\footnote { We have introduced this
notation for convenience.  For example, $h\in \mathcal{V}
\Rightarrow X\subset\subset (0,1) $ means that there exists a
compact set $A \subset (0,1) $ such that $h\in \mathcal{V}
\Rightarrow X\in A $, and similarly for the other variables.}  For
the remainder of this discussion we will restrict our attention to
the dynamics of the system while the slow variables remain in the
set $\mathcal{V} $. To this end, we define the stopping time
$T_\varepsilon(z) =T_\varepsilon =\inf \{\tau\geq 0: \bar {h}
(\tau)\notin \mathcal{V} \text { or } h_\varepsilon(\tau
/\varepsilon) \notin \mathcal{V} \}$.

\begin {thm}
\label{thm:1Dpiston1}

For each $T>0$,
\[
\sup_{\substack{\text{initial conditions}\\
\text{s.t. } h_\varepsilon(0)\in  \mathcal{V} }}\;
\sup_{0\leq\tau\leq T\wedge
T_\varepsilon}\abs{h_\varepsilon(\tau/\varepsilon)-\bar{h}(\tau)}
=\mathcal{O} (\varepsilon  ) \text { as } \varepsilon=M^
{-1/2}\rightarrow 0.
\]
\end{thm}

Note that the stopping time does not unduly restrict the result.
Given any $c $ such that $h=c\Rightarrow X\in (0,1),\: s_{i,j}\in
(0,\infty)$, then by an appropriate choice of the compact set
$\mathcal{V} $ we may ensure that, for all $\varepsilon $
sufficiently small and all initial conditions in our phase space
with $h_\varepsilon(0) =c$, $T_\varepsilon \geq T $. We do this by
choosing $\mathcal{V}\ni c $ such that the distance between
$\partial \mathcal{V} $ and the periodic orbit $\bar h (\tau) $ with
$\bar h (0) =c $ is positive. Call this distance $d $. Then
$T_\varepsilon $ can only occur before $T$ if $h_\varepsilon
(\tau/\varepsilon) $ has deviated by at least $d $ from $ \bar h
(\tau) $ for some $\tau\in [0,T) $. Since the size of the deviations
tends to zero uniformly with $\varepsilon$, this is impossible for
all small $\varepsilon $.

\subsection{The soft core piston problem}
\label {sct:soft_core_results}

In this section, we consider the same system of one piston and gas
particles inside the unit interval considered in Section
\ref{sct:hard_core_results}, but now the interactions of the gas
particles with the walls and with the piston are smooth.  Let
$\kappa\colon\mathbb {R}\rightarrow\mathbb {R} $ be a $\mathcal{C}
^2$ function satisfying
\begin {itemize}
    \item
         $\kappa (x)= 0 $ if $x\geq 1 $,
    \item
        $\kappa ' (x) < 0 $ if $x < 1 $.
\end {itemize}
Let $\delta >0 $ be a parameter of smoothing, and set $\kappa
_\delta (x) =\kappa (x/\delta) $.  Then consider the Hamiltonian
system obtained by having the gas particles interact with the piston
and the walls via the potential
\[
    \sum_{j=1}^{n_1}\kappa_\delta(x_{1,j})+\kappa_\delta(X-x_{1,j})+
    \sum_{j=1}^{n_2}\kappa_\delta(x_{2,j}-X)+\kappa_\delta(1-x_{2,j}).
\]
As before, we set $\varepsilon =M^ {-1/2} $ and $W =V/\varepsilon$.
If we let
\begin {equation}
\label{eq:soft_energies}
\begin {split}
    E_{1,j}&=\frac{1}{2}m_{1,j}v_{1,j}^2+\kappa_\delta(x_{1,j})
        +\kappa_\delta(X-x_{1,j}),\qquad 1\leq j\leq n_1,\\
    E_{2,j} & =\frac{1}{2}m_{2,j}v_{2,j}^2+\kappa_\delta(x_{2,j}-X)
        +\kappa_\delta(1-x_{2,j}),\qquad 1\leq j\leq n_2,\\
\end {split}
\end {equation}
then $E_{i,j}$ may be thought of as the energy associated with a gas
particle, and $W ^2/2+\sum_{j=1}^{n_1}E_{1,j}
+\sum_{j=1}^{n_2}E_{2,j} $ is the conserved energy.

When $\varepsilon =0 $, the Hamiltonian system admits $n_1+n_2+2$
independent first integrals, which we choose this time as $h = (X,
W,E_{1,j},E_{2,j}) $.  While discussing the soft core dynamics we
use the energies $E_{i,j} $ rather than the variables $s_{i,j}  =
\sqrt {2E_{i,j} /m_{i,j} } $, which we used for the hard core
dynamics, for convenience.

For comparison with the hard core results, we formally consider the
dynamics described by setting $\delta=0$ to be the hard core
dynamics described in Section \ref{sct:hard_core_results}. This is
reasonable because we will only consider gas particle energies below
the barrier height $\kappa(0) $.  Then for any
$\varepsilon,\delta\geq 0$, $h_\varepsilon^\delta (t) $ denotes the
actual time evolution of the slow variables. While discussing the
soft core dynamics we often use $\delta $ as a superscript to
specify the dynamics for a certain value of $\delta $. We usually
suppress the dependence on $\delta $, unless it is needed for
clarity.

\subsubsection{Main result in the soft core setting}
\label {sct:1d_sm_results}

We have already seen that when $\delta=0$, there is an appropriate
averaged vector field $\bar H^0$ whose solutions approximate the
actual motions of the slow variables, $h_\varepsilon^0 (t) $. We
will show that when $\delta >0$, there is also an appropriate
averaged vector field $\bar H^\delta$ whose solutions still
approximate the actual motions of the slow variables,
$h_\varepsilon^\delta (t) $. We delay the derivation of $\bar
H^\delta$ until Section \ref{sct:smooth_1D_derivation}.

Fix a compact set $\mathcal{V}\subset\mathbb{R}^{n_1+n_2+2}$ such
that $h\in\mathcal{V}\Rightarrow X\subset\subset (0,1),
W\subset\subset\mathbb{R}$, and $E_{i,j}\subset\subset (0,\kappa
(0)) $ for each $i$ and $j$.  For each $\varepsilon,\delta\geq 0$ we
define the functions $\bar h ^\delta (\cdot) $ and
$T_\varepsilon^\delta $ on our phase space by letting $\bar h
^\delta (\tau) $ be the solution of
\begin {equation}
\label {eq:smooth_1D_averaged_eq}
    \frac {d\bar{h} ^\delta}{d\tau}
    =\bar {H} ^\delta (\bar {h} ^\delta),\qquad \bar {h} ^\delta (0)
    =h_\varepsilon ^\delta (0),
\end {equation}
and $T_\varepsilon ^\delta =\inf \{\tau\geq 0: \bar {h} ^\delta
(\tau)\notin \mathcal{V} \text { or } h_\varepsilon ^\delta (\tau
/\varepsilon) \notin \mathcal{V} \}$.

\begin {thm}
\label{thm:1D_smooth_uniform}

There exists $\delta_0 >0$ such that the averaged vector field $\bar
H^\delta (h) $ is $\mathcal{C} ^1$ on the domain
$\{(\delta,h):0\leq\delta \leq\delta_0,h\in\mathcal{V}\}$.
Furthermore, for each $T>0$,
\[
\sup_{0\leq\delta\leq\delta_0}\; \sup_{\substack {\text {initial conditions}\\
\text {s.t. } h_\varepsilon^\delta(0)\in  \mathcal{V} }}\;
\sup_{0\leq\tau\leq T\wedge T_\varepsilon ^\delta}\abs{h_\varepsilon
^\delta (\tau/\varepsilon)-\bar{h} ^\delta (\tau)} =\mathcal{O}
(\varepsilon  )\text { as }\varepsilon=M^ {-1/2}\rightarrow 0.
\]
\end{thm}

As in Section \ref{sct:1dres}, for any fixed $c $ there exists a
suitable choice of the compact set $\mathcal{V} $ such that for all
sufficiently small $\varepsilon$ and $\delta$, $T_\varepsilon
^\delta\geq T $ whenever $h_\varepsilon ^\delta (0) =c$.

\subsection{Applications and generalizations}

\subsubsection{Relationship between the hard core and the soft core
piston} \label {sct:relationship}

It is not \emph{a priori} clear that we can compare the motions of
the slow variables on the time scale $1/\varepsilon$ for $\delta>0$
versus $\delta=0$, i.e.~compare the motions of the soft core piston
with the motions of the hard core piston on a relatively long time
scale. It is impossible to compare the motions of the fast-moving
gas particles on this time scale as $\varepsilon\rightarrow 0$.  As
we see in Section \ref{sct:sc_proof}, the frequency with which a gas
particle hits the piston changes by an amount $\mathcal{O} (\delta)
$ when we smooth the interaction. Thus, on the time scale
$1/\varepsilon$, the number of collisions is altered by roughly
$\mathcal{O} (\delta/\varepsilon) $, and this number diverges if
$\delta$ is held fixed while $\varepsilon\rightarrow 0$.

Similarly, one might expect that it is impossible to compare the
motions of the soft and hard core pistons as $\varepsilon\rightarrow
0$ without letting $\delta\rightarrow 0$ with $\varepsilon$.
However, from Gronwall's Inequality it follows that if $\bar
h^\delta (0) =\bar h^0(0) $, then $\sup_{0\leq\tau\leq T\wedge
T_\varepsilon ^\delta\wedge T_\varepsilon ^0} \abs{\bar h^\delta
(\tau) -\bar h^0(\tau)}=\mathcal{O} (\delta) $.  From the triangle
inequality and Theorems \ref{thm:1Dpiston1} and
\ref{thm:1D_smooth_uniform} we obtain the following corollary, which
allows us to compare the motions of the hard core and the soft core
piston.
\begin {cor}
\label {cor: 1d_smooth_comparison}

As $\varepsilon=M^ {-1/2},\delta \rightarrow 0$,
\[
\sup_{c\in \mathcal{V} }\; \sup_{\substack {\text {initial conditions}\\
\text {s.t. } h_\varepsilon^\delta(0)=c=h_\varepsilon^0(0) }}\;
 \sup_{0\leq t\leq (T\wedge T_\varepsilon
^\delta\wedge T_\varepsilon ^0)/\varepsilon}\abs{h_\varepsilon
^\delta (t)-h_\varepsilon ^0 (t)} =\mathcal{O} (\varepsilon
)+\mathcal{O} (\delta).
\]
\end {cor}

This shows that, provided the slow variables have the same initial
conditions,
\[
\sup_{0\leq t\leq 1/\varepsilon}\abs{h_\varepsilon ^\delta
(t)-h_\varepsilon ^0 (t)} =\mathcal{O} (\varepsilon )+\mathcal{O}
(\delta).
\]
Thus the motions of the slow variables converge on the time scale
$1/\varepsilon $ as $\varepsilon,\delta\rightarrow 0 $, and it is
immaterial in which order we let these parameters tend to zero.

\subsubsection{The adiabatic piston problem}

We comment on what Theorem \ref{thm:1Dpiston1} says about the
adiabatic piston problem. The initial conditions of the adiabatic
piston problem require that $W(0) =0$. Although our system is so
simple that a proper thermodynamical pressure is not defined, we can
define the pressure of a gas to be the average force received from
the gas particles by the piston when it is held fixed,
i.e.~$P_1=\sum_{j=1} ^{n_1}2m_{1,j} s_{1,j} \frac{s_{1,j}}{2X} =2
E_1/X $ and $P_2=2E_2/(1-X) $. Then if $P_1(0) >P_2(0) $, the
initial condition for our averaged equation \eqref{eq:1davg} has the
motion of the piston starting at the left turning point of a
periodic orbit determined by the effective potential well.  Up to
errors not much bigger than $M^ {-1/2} $, we see the piston
oscillate periodically on the time scale $M^ {1/2}$. If $P_1(0) <
P_2(0) $, the motion of the piston starts at a right turning point.
However, if $P_1(0) = P_2(0) $, then the motion of the piston starts
at the bottom of the effective potential well. In this case of
mechanical equilibrium, $\bar h(\tau) =\bar h(0) $, and we conclude
that, up to errors not much bigger than $M^ {-1/2} $, we see no
motion of the piston on the time scale $M^ {1/2}$.  A much longer
time scale is required to see if the temperatures equilibrate.

\subsubsection{Generalizations}

There are several other ways this work could be generalized. For
example, one could replace the walls at $0$ and $1$ by heat baths,
so that whenever a gas particle reaches a heat bath, it is absorbed,
and another gas particle is emitted with its velocity chosen
independently and randomly according to some distribution.
Similarly, the walls could be replaced by a forcing mechanism (much
like the bumpers in a pinball machine) that changes a colliding gas
particle's kinetic energy.

A simple generalization of Theorem \ref{thm:1Dpiston1}, proved by
similar techniques, follows. The system consists of $N-1$ pistons,
that is, heavy point particles, located inside the unit interval at
positions $X_1< X_2<\dotsc< X_{N-1} $. Walls are located at
$X_0\equiv 0$ and $X_N\equiv 1$, and the piston at position $X_i$
has mass $M_i$. Then the pistons divide the unit interval into $N$
chambers. Inside the $i^ {th}$ chamber, there are $n_i\geq 1$ gas
particles whose locations and masses will be denoted by $x_{i,j}$
and $m_{i,j}$, respectively, where $1\leq j\leq n_i$. All of the
particles are point particles, and the gas particles interact with
the pistons and with the walls via elastic collisions. However, the
gas particles do not directly interact with each other. We scale the
piston masses as $M_i=\hat M_i/\varepsilon^2$ with $\hat M_i$
constant, define $W_i$ by $dX_i/dt=\varepsilon W_i$, and let $E_i$
be the kinetic energy of the gas particles in the $i^ {th}$ chamber.
Then we can find an appropriate averaged equation whose solutions
have the pistons moving like an $N$-dimensional particle inside a
potential well with an effective Hamiltonian
\[
    \frac {1} {2}\sum_{i=1}^{N-1}\hat M_i W_i^2 +
    \sum_{i=1}^{N}
    \frac{E_i(0)(X_i(0)-X_{i-1}(0))^2}{(X_i-X_{i-1})^2}.
\]
If we write the slow variables as $h= (X_i,W_i,\abs{v_{i,j}})$ and
fix a compact set $\mathcal{V}$ such that $h\in \mathcal{V}
\Rightarrow X_{i+1} -X_i\subset\subset (0,1),W_i\subset\subset
\mathbb R$, and $\abs{v_{i,j}} \subset\subset (0,\infty)$, then the
convergence of the actual motions of the slow variables to the
averaged solutions is exactly the same as the convergence given in
Theorem \ref{thm:1Dpiston1}.

\section{Preparatory material}

Before the proofs of our main results, we present some averaging
results, as well as a heuristic derivation of the averaged equation
for the hard core piston.  This material provides background for our
work and establishes some notation.

\subsection{The averaging framework}
\label {sct:anos}

In this section, consider a family of ordinary differential
equations
\begin {equation}\label{eq:ode}
\frac {dz} {dt}=Z(z,\varepsilon)
\end {equation}
on a smooth, finite-dimensional Riemannian manifold $\mathcal{M}$,
which depends on the real parameter $\varepsilon\in
[0,\varepsilon_0]$. Assume
\begin {itemize}
\item
     {\em  Regularity:} the functions $Z$ and
            $\partial Z/\partial \varepsilon$ are both
            $\mathcal{C}^1$ on $\mathcal{M}\times [0,\varepsilon_0] $.
\end {itemize}
We denote the flow generated by $Z(\cdot,\varepsilon) $ by
$z_\varepsilon(t,z)=z_\varepsilon(t) $. We will usually suppress the
dependence on the initial condition $z=z_\varepsilon(0,z) $.  Think
of $z_\varepsilon(\cdot) $ as being a random variable whose domain
is the space of initial conditions for the differential equation
\eqref{eq:ode} and whose range is the space of continuous paths
(depending on the parameter $t$) in $\mathcal{M}$.
\begin {itemize}
\item
     {\em  Existence of smooth integrals:} $z_0(t) $ has $m$ independent
            $\mathcal{C}^2$ first integrals $h=(h_1,\dotsc,h_m)$.
\end {itemize}
Then $h$ is conserved by $z_0(t) $, and at every point the linear
operator $\partial h/\partial z$ has full rank.  It follows from the
implicit function theorem that each level set $\mathcal{M}_c=
\{h=c\} $ is a smooth submanifold of co-dimension $m$, which is
invariant under $z_0(t) $. Further, assume that there exists an open
ball $\mathcal{U}\subset\mathbb{R}^m$ satisfying:
\begin {itemize}
    \item
        {\em Compactness:} $\forall c\in\mathcal{U},\: \mathcal{M}_c$ is compact.
    \item
        {\em Preservation of smooth measures:} $\forall c\in\mathcal{U}$,
        $z_0(t)\arrowvert_{\mathcal{M}_c} $ preserves a smooth measure
        $\mu_c$ that varies smoothly with $c$, i.e.~there exists a
        $\mathcal{C}^1 $ function $f:\mathcal{M}\rightarrow\mathbb{R}_{>0} $ such that
        $f\arrowvert_{\mathcal{M}_c} $ is the density of $\mu_c$ with respect
        to the restriction of Riemannian volume.
\end {itemize}

Set $h_\varepsilon(t) =h(z_\varepsilon(t)) $.  Since $dh_0/dt
\equiv0 $, Hadamard's Lemma allows us to write
\[
\frac {dh_\varepsilon} {dt} =\varepsilon
H(z_\varepsilon,\varepsilon)
\]
for some function $H:\mathcal{M}\times
[0,\varepsilon_0]\rightarrow\mathcal{U}$. Observe that
\[
    \frac {dh_\varepsilon} {dt}(t) =
    Dh(z_\varepsilon(t))Z(z_\varepsilon(t),\varepsilon)=
    Dh(z_\varepsilon(t))\bigl(Z(z_\varepsilon(t),\varepsilon)
    -Z(z_\varepsilon(t),0)\bigr),
\]
so that
\[
H(z,0) =\mathcal{L}_{\frac {\partial Z}
{\partial\varepsilon}\arrowvert_{\varepsilon=0}}h.
\]
Here $\mathcal{L}$ denotes the Lie derivative.

Define the averaged vector field $\bar{H} $ by
\begin{equation}
\label{eq:Anosov_ avg}
    \bar {H} (h)
    =\int_{\mathcal{M}_h}H(z,0)d\mu_h(z).
\end{equation}
Then $\bar {H} $ is $\mathcal{C}^1$. Fix a compact set
$\mathcal{V}\subset\mathcal{U}$, and let $\bar{h} (\tau) $ be the
solution of
\[
\frac {d\bar{h}}{d\tau} =\bar {H} (\bar {h}),\qquad \bar {h} (0)
=h_\varepsilon(0).
\]
We only consider the dynamics in a compact subset of phase space, so
for initial conditions $z\in h^ {-1} \mathcal{U}$, define the
stopping time $T_\varepsilon=\inf \{\tau\geq 0: \bar {h}
(\tau)\notin \mathcal{V} \text { or } h_\varepsilon(\tau
/\varepsilon) \notin \mathcal{V} \}$.

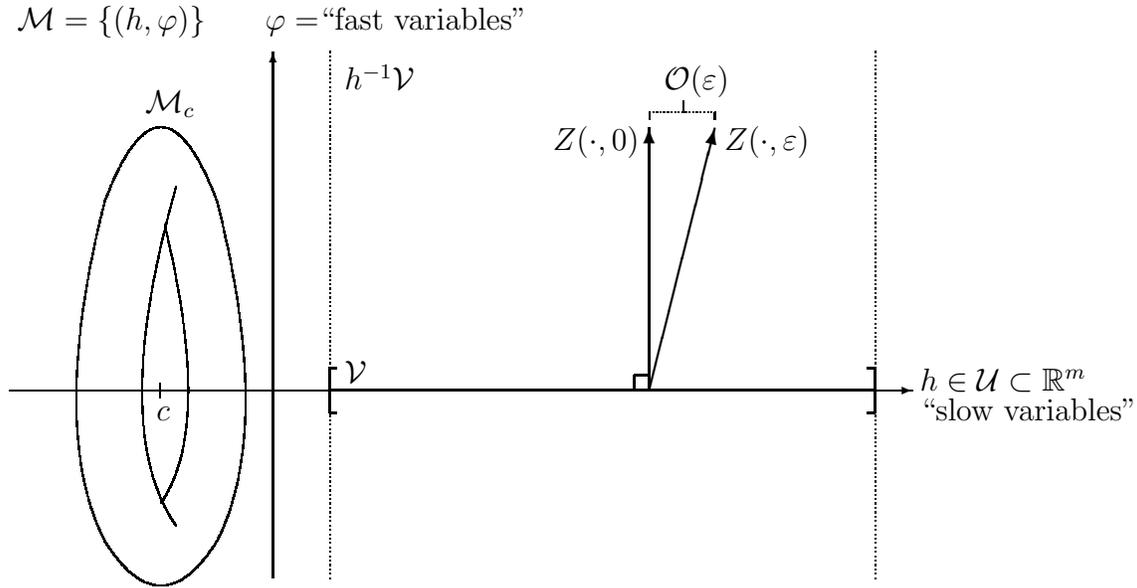
\begin{figure}
    \begin {center}
    \setlength{\unitlength}{1 cm}
    \begin{picture}(15,8)
        \put(0.1,7.3){$\mathcal{M}=\{(h,\varphi)\}$}
        \put(0,2.5){\vector(1,0){12}}
        \put(12.1,2.5){$h\in\mathcal{U}\subset\mathbb{R}^m$}
        \put(12.1,2.1){``slow variables''}
        \put(3.5,0){\vector(0,1){7}}
        \put(3.4,7.3){$\varphi=$``fast variables''}
        \thicklines
        \put(8.5,2.5){\vector(0,1){3.5}}
        \put(8.5,2.5){\vector(1,4){0.87}}
        \put(8.3,2.7){\line(1,0){.2}} \put(8.3,2.5){\line(0,1){.2}}
        \put(7.22,5.7){$Z(\cdot,0) $}
        \put(9.5,5.7){$Z(\cdot,\varepsilon) $}
        \thinlines
        \qbezier[18](8.5,6.2)(8.93,6.2)(9.37,6.2)
        \put(8.95,6.2){\line(0,1){.15}}
        \put(8.5,6.2){\line(0,-1){0.1}}
        \put(9.37,6.2){\line(0,-1){0.1}}
        \put(8.7,6.5){$\mathcal{O} (\varepsilon)$}
        \put(2,2.6){\line(0,-1){0.2}}
        \put(1.95,2.1){$c$}
        \qbezier(2.2,5.2)(1.3,2)(2.2,0.7)
        \qbezier(2.05,4.7)(2.7,2)(2,1)
        \qbezier(1.25,5)(.5,2)(1.25,0.5)
        \qbezier(2.75,5)(3.5,2)(2.75,0.5)
        \qbezier(1.25,5)(2,7)(2.75,5)
        \qbezier(1.25,0.5)(2,-0.7)(2.75,0.5)
        \put(1.8,6.2){$\mathcal{M}_c$}
        \qbezier[120](4.25,0)(4.25,3.5)(4.25,7)
        \qbezier[120](11.5,0)(11.5,3.5)(11.5,7)
        \put(4.45,2.6){$\mathcal{V}$}
        \put(4.45,6.5){$h^ {-1}\mathcal{V}$}
        \thicklines
        \put(4.25,2.5){\line(1,0){7.25}}
        \put(4.25,2.2){\line(0,1){.6}}
        \put(4.25,2.2){\line(1,0){.1}}
        \put(4.25,2.8){\line(1,0){.1}}
        \put(11.5,2.2){\line(0,1){.6}}
        \put(11.5,2.2){\line(-1,0){.1}}
        \put(11.5,2.8){\line(-1,0){.1}}
    \end{picture}
    \end {center}
    \caption{A schematic of the phase space $\mathcal{M}$.}\label{fig:phase_space}
\end{figure}

Heuristically, think of the phase space $\mathcal{M}$ as being a
fiber bundle whose base is the open set $\mathcal{U}$ and whose
fibers are the compact sets $\mathcal{M}_h$. See Figure
\ref{fig:phase_space}. Then the vector field $Z(\cdot,0) $ is
perpendicular to the base, so its orbits $z_0(t) $ flow only along
the fibers.  Now when $0< \varepsilon\ll 1 $, the vector field
$Z(\cdot,\varepsilon) $ acquires a component of size $\mathcal{O}
(\varepsilon) $ along the base, and so its orbits $z_\varepsilon(t)
$ have a small drift along the base, which we can follow by
observing the evolution of $h_\varepsilon(t) $. Because of this, we
refer to $h$ as consisting of the slow variables. Other variables,
used to complete $h$ to a parameterization of (a piece of) phase
space, are called fast variables. Note that $h_\varepsilon(t)$
depends on all the dimensions of phase space, and so it is not the
flow of a vector field on the $m$-dimensional space $\mathcal{U}$.
However, because the motion along each fiber is relatively fast
compared to the motion across fibers, we hope to be able to average
over the fast motions and obtain a vector field on $\mathcal{U}$
that gives a good description of $h_\varepsilon(t)$ over a
relatively long time interval, independent of where the solution
$z_\varepsilon(t) $ started on $\mathcal{M}_{h_\varepsilon(0)} $.
Because our averaged vector field, as defined by Equation
\eqref{eq:Anosov_ avg}, only accounts for deviations of size
$\mathcal{O}(\varepsilon) $, we cannot expect this time interval to
be longer than size $\mathcal{O}(1/\varepsilon)$.  In terms of the
slow time $\tau=\varepsilon t$, this length becomes $\mathcal{O}(1)
$.  In other words, the goal of the first-order averaging method
described above should be to show that, in some sense,
$\sup_{0\leq\tau\leq 1\wedge
T_\varepsilon}\abs{h_\varepsilon(\tau/\varepsilon)-\bar{h}(\tau)}\rightarrow
0$ as $\varepsilon\rightarrow 0$.  This is often referred to as the
averaging principle.

\subsubsection{Some averaging results}

So far, we are in a general averaging setting.  Frequently, one also
assumes that the invariant submanifolds, $\mathcal{M}_h$, are tori,
and that there exists a choice of coordinates $z=(h,\varphi) $ on
$\mathcal{M} $ in which the differential equation \eqref{eq:ode}
takes the form
\[
    \frac {dh} {dt}=\varepsilon H(h,\varphi,\varepsilon),\qquad
    \frac {d\varphi} {dt}=\Phi(h,\varphi,\varepsilon).
\]
Then if $\varphi\in S^1$ and the differential equation for the fast
variable is regular, i.e.~$\Phi(h,\varphi,0)$ is bounded away from
zero for $h\in\mathcal{U} $,
\[
    \sup_{\substack {\text {initial conditions}\\
    \text {s.t. } h_\varepsilon(0)\in  \mathcal{V} }}\;
    \sup_{0\leq\tau\leq 1\wedge
    T_\varepsilon}\abs{h_\varepsilon(\tau/\varepsilon)-\bar{h}(\tau)}
    =\mathcal{O} (\varepsilon )\text { as }\varepsilon\rightarrow 0.
\]
See for example Chapter 5 in \cite{SV85} or Chapter 3 in
\cite{LM88}.

When the differential equation for the fast variable is not regular,
or when there is more than one fast variable, the typical averaging
result becomes much weaker than the uniform convergence above.  For
example, consider the case when $\varphi\in\mathbb{T}^n$, $n>1$, and
the unperturbed motion is quasi-periodic,
i.e.~$\Phi(h,\varphi,0)=\Omega(h) $.  Also assume that
$H\in\mathcal{C} ^ {n+2}$ and that $\Omega$ is nonvanishing and
satisfies a nondegeneracy condition on $\mathcal{U} $ (for example,
$\Omega :\mathcal{U}\rightarrow\mathbb{T}^n$ is a submersion).  Let
$P$ denote Riemannian volume on $\mathcal{M} $.
Neishtadt~\cite{LM88,Nei76} showed that in this situation, for each
fixed $\delta>0$,
\[
    P\left (\sup_{0\leq\tau\leq 1\wedge
    T_\varepsilon}\abs{h_\varepsilon(\tau/\varepsilon)-\bar{h}(\tau)}
    \geq \delta
    \right)
    =\mathcal{O} (\sqrt\varepsilon /\delta),
\]
and that this result is optimal.  Thus, the averaged equation only
describes the actual motions of the slow variables in probability on
the time scale $1/\varepsilon $ as $\varepsilon\rightarrow 0$.


Neishtadt's result was motivated by a general averaging theorem for
smooth systems due to Anosov. This theorem requires none of the
additional assumptions in the averaging results above. Under the
conditions of regularity, existence of smooth integrals,
compactness, and preservation of smooth measures, as well as
\begin {itemize}
    \item
        {\em Ergodicity:} for Lebesgue almost every
        $c\in\mathcal{U}$, $(z_0(\cdot),\mu_c) $ is ergodic,
\end {itemize}
Anosov~\cite{Ano60,LM88} showed that
\[
\sup_{0\leq\tau\leq 1\wedge
T_\varepsilon}\abs{h_\varepsilon(\tau/\varepsilon)-\bar{h}(\tau)}\rightarrow
0\text { in probability (w.r.t. Riemannian volume) as
}\varepsilon\rightarrow 0.
\]

If we consider $h_\varepsilon(\cdot) $ and $\bar h(\cdot) $ to be
random variables, this is a version of the weak law of large
numbers. In general, we can do no better: There is no general strong
law in this setting. There exists an example due to Neishtadt (which
comes from the equations for the motion of a pendulum with linear
drag being driven by a constant torque) where for no initial
condition in a positive measure set do we have convergence of
$h_\varepsilon(t)$ to $\bar{h}(\varepsilon t)$ on the time scale
$1/\varepsilon $ as $\varepsilon\rightarrow 0$~\cite{Kif04}. Here,
the phase space is $\mathbb{R}\times S^1$, and the unperturbed
motion is (uniquely) ergodic on all but one fiber.

In 2004, Kifer~\cite{Kif04b} gave necessary and sufficient
conditions for the averaging principle to hold in an averaged sense
with respect to initial conditions.  He also showed explicitly that
his conditions are met in the setting of Anosov's theorem.  A
recent, relatively simple proof of the version of Anosov's theorem
stated above is due to Dolgopyat~\cite {Dol05}.

\subsection{Heuristic derivation of the averaged equation for the
        hard core piston}
\label {sct:1dps_hc_heuristic}

We present here a heuristic derivation of Sinai's averaged equation
\eqref{eq:1davg} that is found in~\cite {Dol05}.

First, we examine interparticle collisions when $\varepsilon>0$.
When a particle on the left, say the one at position $x_{1,j} $,
collides with the piston, $s_{1,j} $ and $W$ instantaneously change
according to the laws of elastic collisions:
\begin {equation}\label {eq:v1lin}
    \begin{bmatrix}
    v_{1,j}^+\\ V^+
    \end{bmatrix}
    =
    \frac{1}{m_{1,j}+M}
    \begin{bmatrix}
    m_{1,j}-M& 2M\\
    2m_{1,j}& M-m_{1,j}\\
    \end{bmatrix}
    \begin{bmatrix}
    v_{1,j}^-\\ V^-
    \end{bmatrix}.
\end {equation}
If the speed of the left gas particle is bounded away from zero, and
$W=M^ {1/2} V $ is also bounded, it follows that for all
$\varepsilon$ sufficiently small, any collision will have $v_{1,j}^
-
>0$ and $v_{1,j}^ + <0$. In this case, when we translate Equation
\eqref{eq:v1lin} into our new coordinates, we find that
\begin {equation}\label {eq:s1lin}
    \begin{bmatrix}
    s_{1,j}^+\\ W^+
    \end{bmatrix}
    =
    \frac{1}{1+\varepsilon^2 m_{1,j}}
    \begin{bmatrix}
    1-\varepsilon^2 m_{1,j}& - 2\varepsilon\\
    2\varepsilon m_{1,j}& 1-\varepsilon^2 m_{1,j}\\
    \end{bmatrix}
    \begin{bmatrix}
    s_{1,j}^-\\ W^-
    \end{bmatrix},
\end {equation}
so that
\[
\begin {split}
    \Delta s_{1,j}&=s_{1,j}^+ -s_{1,j}^- =-2\varepsilon W^-
    +\mathcal{O}(\varepsilon^2),\\
    \Delta W &=W^+ -W^- =+2\varepsilon m_{1,j}
    s_{1,j}^- +\mathcal{O}(\varepsilon^2).\\
\end {split}
\]

The situation is analogous when particles on the right collide with
the piston. For all $\varepsilon $ sufficiently small, $s_{2,j} $
and $W$ instantaneously change by
\[
\begin {split}
    \Delta W &=W^+ -W^- = - 2\varepsilon m_{2,j}
    s_{2,j}^- +\mathcal{O}(\varepsilon^2),\\
    \Delta s_{2,j}&=s_{2,j}^+ -s_{2,j}^- = + 2\varepsilon W^-
    +\mathcal{O}(\varepsilon^2).\\
\end {split}
\]
We defer discussing the rare events in which multiple gas particles
collide with the piston simultaneously, although we will see that
they can be handled appropriately.

Let $\Delta t $ be a length of time long enough such that the piston
experiences many collisions with the gas particles, but short enough
such that the slow variables change very little, in this time
interval.  From each collision with the particle at position
$x_{1,j} $, $W$ changes by an amount $+2\varepsilon m_{1,j}
s_{1,j}+\mathcal{O}(\varepsilon^2)$, and the frequency of these
collisions is approximately $\frac {s_{1,j}} {2X} $. Arguing
similarly for collisions with the other particles, we guess that
\[
    \frac {\Delta W} {\Delta t} =
    \varepsilon\sum_{j=1} ^ {n_1} 2m_{1,j}
    s_{1,j}\frac{s_{1,j}}{2X}
    - \varepsilon\sum_{j=1} ^ {n_2}
    2m_{2,j} s_{2,j}\frac{s_{2,j}}{2(1-X)}
    +\mathcal{O}(\varepsilon^2).
\]
Note that not only does the position of the piston change slowly in
time, but its velocity also changes slowly, i.e. the piston has
inertia. With $\tau=\varepsilon t$ as the slow time, a reasonable
guess for the averaged equation for $W$ is
\[
    \frac {dW}{d\tau}=\frac{\sum_{j=1} ^ {n_1}m_{1,j} s_{1,j}^2}{X}
    -\frac{\sum_{j=1} ^ {n_2}m_{2,j} s_{2,j}^2}{1-X}.
\]
Similar arguments for the other slow variables lead to the averaged
equation \eqref{eq:1davg}.

\section{Proof of the main result for the hard core piston}
\label {sct:hc_proof}

\subsection{Proof of Theorem \ref{thm:1Dpiston1} with only one gas
            particle on each side}

We specialize to the case when there is only one gas particle on
either side of the piston, i.e.~we assume that $ n_1= n_2=1$.  We
then denote $x_{1,1} $ by $x_1$, $m_{2,2} $ by $m_2$, etc. This
allows the proof's major ideas to be clearly expressed, without
substantially limiting their applicability. At the end of this
section, we outline the simple generalizations needed to make the
proof apply in the general case.

\subsubsection{A choice of coordinates on the phase space for a
    three particle system}
\label {sct:1dps_hcPS}

As part of our proof, we choose a set of coordinates on our
six-dimensional phase space such that, in these coordinates, the
$\varepsilon=0$ dynamics are smooth.  Complete the slow variables
$h= (X,W,s_1,s_2)$ to a full set of coordinates by adding the
coordinates $\varphi _i\in [0,1]/ \, 0\sim 1=S^1,\:i=1,2$, defined
as follows:
\[
\begin {split}
    \varphi_1 =\varphi_1 (x_1,v_1,X)=&\begin {cases}
    \frac{x_1}{2X} &\text { if } v_1 >0\\
    1-\frac{x_1}{2X} &\text { if } v_1<0\\
    \end {cases}
    \\
    \varphi_2=\varphi_2(x_2,v_2,X)=&\begin {cases}
    \frac{1-x_2}{2(1-X)} &\text { if } v_2<0\\
    1- \frac{1-x_2}{2(1-X)} &\text { if } v_2>0\\
    \end {cases}\\
\end {split}.
\]
When $\varepsilon =0 $, these coordinates are simply the angle
variable portion of action-angle coordinates for an integrable
Hamiltonian system.  They are defined such that collisions occur
between the piston and the gas particles precisely when $\varphi_1 $
or $\varphi_2 =1/2$. Then $z = (h,\varphi_1,\varphi_2) $ represents
a choice of coordinates on our phase space, which is homeomorphic to
$(\text{a subset of } \mathbb{R}^4)\times \mathbb{ T} ^ 2 $.  We
abuse notation and also let $h(z) $ represent the projection onto
the first four coordinates of $z$.

Now we describe the dynamics of our system in these coordinates.
When $\varphi_1,\varphi_2\ne 1/2$,
\[
\begin {split}
    \frac{d\varphi_1}{dt} =&\begin {cases}
    \frac{s_1}{2X}-\frac{\varepsilon W}{X}\varphi_1&\text { if }
    0\leq \varphi_1<1/2\\
    \frac{s_1}{2X}+\frac{\varepsilon W}{X}(1-\varphi_1)&\text { if }
    1/2< \varphi_1\leq 1\\
    \end {cases}
    \\
    \frac{d\varphi_2}{dt} =&\begin {cases}
    \frac{s_2}{2(1-X)}+\frac{\varepsilon W}{1-X}\varphi_2&\text { if }
    0\leq \varphi_2<1/2\\
    \frac{s_2}{2(1-X)}-\frac{\varepsilon W}{1-X}(1-\varphi_2)&\text { if }
    1/2< \varphi_2\leq 1\\
    \end {cases}
\end {split}.
\]
Hence between interparticle collisions, the dynamics are smooth and
are described by
\begin {equation}\label {eq:1dode}
\begin {split}
    \frac{dX}{dt}&=\varepsilon W,\\
    \frac{dW}{dt}&=
    \frac{ds_1}{dt}=
    \frac{ds_2}{dt}=0,\\
    \frac{d\varphi_1}{dt}&=\frac{s_1}{2X}+
        \mathcal{O}(\varepsilon),\\
    \frac{d\varphi_2}{dt}&=\frac{s_2}{2 (1 - X)}+
        \mathcal{O}(\varepsilon).\\
\end {split}
\end {equation}

When $\varphi_1$ reaches $1/2$, while $\varphi_2\neq 1/2$, the
coordinates $X,s_2,\varphi_1$, and $\varphi_2$ are instantaneously
unchanged, while $s_1$ and $W$ instantaneously jump, as described by
Equation \eqref{eq:s1lin}. It is curious that $s_1^+ +\varepsilon
W^+ = s_1^- -\varepsilon W^-$, so that $d\varphi_1/dt $ is
continuous as $\varphi_1 $ crosses $1/2 $. However, the collision
induces discontinuous jumps of size $\mathcal{O} (\varepsilon ^2) $
in $dX/dt$ and $d\varphi_2/dt$.  Denote the linear transformation in
Equation \eqref{eq:s1lin} with $j=1$ by $A_{1,\varepsilon} $.  Then
$A_{1,\varepsilon} =\begin{bmatrix}
    1& - 2\varepsilon\\
    2\varepsilon m_1& 1\\
    \end{bmatrix} +\mathcal{O} (\varepsilon^2) $.

The situation is analogous when $\varphi_2$ reaches $1/2$, while
$\varphi_1\neq 1/2$.  Then $W$ and $s_2$ are instantaneously
transformed by a linear transformation $A_{2,\varepsilon}
=\begin{bmatrix}
    1&  -2\varepsilon m_2\\
    2\varepsilon& 1\\
    \end{bmatrix} +\mathcal{O} (\varepsilon^2) $.

We also account for the possibility of all three particles colliding
simultaneously.  There is no completely satisfactory way to do this,
as the dynamics have an essential singularity near $\{\varphi_1
=\varphi_2 =1/2\} $.  Furthermore, such three particle collisions
occur with probability zero with respect to the invariant measure
discussed below.  However, the two $3\times 3$ matrices
\[
    \begin{bmatrix}
    A_{1,\varepsilon} & 0\\
    0& 1\\
    \end{bmatrix},
    \begin{bmatrix}
    1 & 0\\
    0& A_{2,\varepsilon}\\
    \end{bmatrix}
\]
have a commutator of size $\mathcal{O} (\varepsilon ^2) $.  We will
see that this small of an error will make no difference to us as
$\varepsilon\rightarrow 0$, and so when $\varphi_1 =\varphi_2 =1/2
$, we pretend that the left particle collides with the piston
instantaneously before the right particle does.  Precisely, we
transform the variables $s_1,\: W, $ and $s_2 $ by
\[
\begin{bmatrix}
    s_1^+\\ W^+ \\ s_2^+\\
    \end{bmatrix}=
\begin{bmatrix}
     1 & 0\\
    0& A_{2,\varepsilon}\\
    \end{bmatrix}
\begin{bmatrix}
    A_{1,\varepsilon} & 0\\
    0& 1\\
\end{bmatrix}
\begin{bmatrix}
    s_1^-\\ W^-\\ s_2^-\\
    \end{bmatrix}.
\]
We find that
\[
\begin {split}
    \Delta s_1&=s_1^+ -s_1^- = - 2\varepsilon W^-
    +\mathcal{O}(\varepsilon^2),\\
    \Delta W &=W^+ -W^- = +2\varepsilon m_1 s_1^-
    - 2\varepsilon m_2
    s_2^- +\mathcal{O}(\varepsilon^2),\\
    \Delta s_2&=s_2^+ -s_2^- = + 2\varepsilon W^-
    +\mathcal{O}(\varepsilon^2).\\
\end {split}
\]

The above rules define a flow on the phase space, which we denote by
$z_\varepsilon (t) $.  We denote its components by
$X_\varepsilon(t),\: W_\varepsilon(t),\: s_{1,\varepsilon} (t), $
etc.  When $\varepsilon
> 0 $, the flow is not continuous, and for
definiteness we take $z_\varepsilon(t) $ to be left continuous in
$t$.

Because our system comes from a Hamiltonian system, it preserves
Liouville
 measure. In our coordinates, this measure has a density proportional
 to
$X (1 -X) $.  That this measure is preserved also follows from the
fact that the ordinary differential equation \eqref{eq:1dode}
preserves this measure, and the matrices
$A_{1,\varepsilon},\:A_{2,\varepsilon} $ have determinant $1 $. Also
note that the set $\{\varphi_1 =\varphi_2 =1/2 \} $ has co-dimension
two, and so $\bigcup_t z_\varepsilon(t)\{\varphi_1 =\varphi_2 =1/2
\} $ has co-dimension one, which shows that only a measure zero set
of initial conditions will give rise to three particle collisions.

\subsubsection{Argument for uniform convergence} \label
{sct:1dps_hc_unif}

\paragraph*{Step 1:  Reduction using Gronwall's Inequality.}

Define $H(z) $ by
\[
    H(z) =
    \begin{bmatrix}
    W\\
    2m_1 s_1 \delta_{\varphi_1=1/2}-2m_2
     s_2\delta_{\varphi_2=1/2}\\
    -2W\delta_{\varphi_1=1/2}\\
    2W\delta_{\varphi_2=1/2}\\
    \end{bmatrix}.
\]
Here we make use of Dirac delta functions.  All integrals involving
these delta functions may be replaced by sums.  We explicitly deal
with any ambiguities arising from collisions occurring at the limits
of integration.

\begin {lem}
For $0 \leq t\leq \frac{T\wedge T_\varepsilon}{\varepsilon}$,
\[
    h_\varepsilon(t)-h_\varepsilon(0)=
    \varepsilon\int_0^t
    H(z_\varepsilon(s))ds+\mathcal{O}(\varepsilon),
\]
where any ambiguity about changes due to collisions occurring
precisely at times $0 $ and $t $ is absorbed in the $\mathcal{O}
(\varepsilon) $ term.
\end {lem}

\begin {proof}
There are four components to verify.  The first component requires
that $X_\varepsilon(t)-X_\varepsilon(0)=\varepsilon\int_0^t
W_\varepsilon(s)ds+\mathcal{O}(\varepsilon)$.  This is trivially
true because $X_\varepsilon(t)-X_\varepsilon(0)=\varepsilon\int_0^t
W_\varepsilon(s)ds$.

The second component states that
\begin {equation}
\label {eq:delt_w}
    W_\varepsilon(t)-
    W_\varepsilon(0)
    =
    \varepsilon\int_0^t 2m_1
    s_{1,\varepsilon}(s) \delta_{\varphi_{1,\varepsilon}(s)=1/2}-2m_2
   s_{2,\varepsilon}(s)\delta_{\varphi_{2,\varepsilon}(s)=1/2}ds
     +\mathcal{O} (\varepsilon).
\end {equation}
Let $r_k $ and $q_j $ be the times in $(0,t) $ such that
$\varphi_{1,\varepsilon}(r_k)=1/2$ and
$\varphi_{2,\varepsilon}(q_j)=1/2$, respectively.  Then
\[
    W_\varepsilon(t)- W_\varepsilon(0)=
     \sum_{r_k}\Delta W_\varepsilon (r_k)+
    \sum_{q_j}\Delta W_\varepsilon (q_j)
     +\mathcal{O} (\varepsilon).
\]
Observe that there exists $\omega >0 $ such that for all
sufficiently small $\varepsilon $ and all $h\in\mathcal{V}$, $
1/\omega < \frac {d\varphi_i} {dt} <\omega $.  Thus the number of
collisions in a time interval grows no faster than linearly in the
length of that time interval.  Because $t\leq T/\varepsilon $, it
follows that
\[
   W_\varepsilon(t)- W_\varepsilon(0)=\\
    \varepsilon\sum_{r_k} 2m_1
    s_{1,\varepsilon}(r_k) -\varepsilon\sum_{q_j}2m_2
   s_{2,\varepsilon}(q_j)
     +\mathcal{O} (\varepsilon),
\]
and Equation \eqref {eq:delt_w} is verified. Note that because
$\mathcal{V} $ is compact, there is uniformity over all initial
conditions in the size of the $\mathcal{O} (\varepsilon) $ terms
above.  The third and fourth components are handled similarly.
\end {proof}

Next, $\bar {h}(\tau) $ satisfies the integral equation
\[
\bar {h}(\tau) -\bar h(0) = \int_0^{\tau}\bar H(\bar
h(\sigma))d\sigma,
\]
while $h_\varepsilon(\tau/\varepsilon)$ satisfies
\[
\begin {split}
    h_\varepsilon(\tau/\varepsilon)-h_\varepsilon(0)
    &=
    \mathcal{O}(\varepsilon) +\varepsilon\int_0^{\tau/\varepsilon}
    H(z_\varepsilon(s))ds\\
    &=\mathcal{O}(\varepsilon) +
    \varepsilon\int_0^{\tau/\varepsilon}
    H(z_\varepsilon(s))-
    \bar H(h_\varepsilon(s))ds+
    \int_0^{\tau}\bar H( h_\varepsilon(\sigma/\varepsilon))d\sigma
\end {split}
\]
for $0\leq\tau\leq T\wedge T_\varepsilon$.

Define $e_\varepsilon(\tau) =\varepsilon\int_0^{\tau/\varepsilon}
H(z_\varepsilon(s))- \bar H(h_\varepsilon(s))ds$. It follows from
Gronwall's Inequality that
\begin {equation}
\label {eq:1dgronwall}
    \sup_{0\leq \tau\leq T\wedge T_\varepsilon}
    \abs{\bar h(\tau)-h_\varepsilon(\tau/\varepsilon)}\leq
    \left(\mathcal{O}(\varepsilon)+
    \sup_{0\leq \tau\leq T\wedge T_\varepsilon}
    \abs{e_\varepsilon(\tau)}\right)e^{ \Lip{\bar
    H\arrowvert _\mathcal{V}}T}.
\end {equation}
Gronwall's Inequality is usually stated for continuous paths, but
the standard proof (found in \cite{SV85}) still works for paths that
are merely integrable, and $\abs{\bar
h(\tau)-h_\varepsilon(\tau/\varepsilon)}$ is piecewise smooth.

\paragraph*{Step 2:  A splitting according to particles.}

Now
\[
    H(z)-\bar H(h) =
    \begin{bmatrix}
    0\\
    2m_1 s_1  \delta_{\varphi_1=1/2}-m_{1} s_{1}^2/X\\
    -2W\delta_{\varphi_1=1/2}+s_1 W/X\\
    0\\
    \end{bmatrix}
    +
    \begin{bmatrix}
    0\\
    -2m_2 s_2\delta_{\varphi_2=1/2}+m_2 s_2^2/(1-X)\\
    0\\
    2W\delta_{\varphi_2=1/2}-s_2 W/(1-X)\\
    \end{bmatrix}
    ,
\]
and so, in order to show that $\sup_{0\leq \tau\leq T\wedge
T_\varepsilon}\abs{e_\varepsilon(\tau)} =\mathcal{O} (\varepsilon)
$, it suffices to show that
\begin {align*}\label{eq:hc_separation}
    &
    \sup_{0\leq \tau\leq
    T\wedge T_\varepsilon}
    \abs{\int_0^{\tau/\varepsilon}
    s_{1,\varepsilon}(s)\delta_{\varphi_{1,\varepsilon}(s)=1/2}-
    \frac{s_{1,\varepsilon}(s)^2}{2X_\varepsilon(s)}ds}
    =\mathcal{O} (1),
    \\
    &
    \sup_{0\leq \tau\leq
    T\wedge T_\varepsilon}
    \abs{\int_0^{\tau/\varepsilon}
    W_{\varepsilon}(s)\delta_{\varphi_{1,\varepsilon}(s)=1/2}-
    \frac{W_{\varepsilon}(s)s_{1,\varepsilon}(s)}{2X_\varepsilon(s)}ds}
    =\mathcal{O} (1),
\end {align*}
as well as two analogous claims about terms involving
$\varphi_{2,\varepsilon} $.  Thus we have effectively separated the
effects of the different gas particles, so that we can deal with
each particle separately. We will only show that $ \sup_{0\leq
\tau\leq
    T\wedge T_\varepsilon}
    \abs{\int_0^{\tau/\varepsilon}
    s_{1,\varepsilon}(s)\delta_{\varphi_{1,\varepsilon}(s)=1/2}-
    \frac{s_{1,\varepsilon}(s)^2}{2X_\varepsilon(s)}ds}
    =\mathcal{O} (1)$.  The other three
terms can be handled similarly.

\paragraph*{Step 3:  A sequence of times adapted for ergodization.}

Ergodization refers to the convergence along an orbit of a
function's time average to its space average.  Because of the
splitting according to particles above, one can easily check that
$\frac {1} {t}\int_{0}^{t} H(z_0(s))ds =\bar H (h_0) +\mathcal{O}
(1/t) $, even when $z_0 (\cdot) $ restricted to the invariant tori
$\mathcal{M}_{h_0} $ is not ergodic.  In this step, for each initial
condition $z_\varepsilon(0) $ in our phase space, we define a
sequence of times $t_{k,\varepsilon} $ inductively as follows:
$t_{0,\varepsilon}=\inf\{t\geq 0:\varphi_{1,\varepsilon}(t)=0\}$,
$t_{k+1,\varepsilon}=\inf\{t>t_{k,\varepsilon}:\varphi_{1,\varepsilon}(t)=0\}$.
This sequence is chosen because $\delta_{\varphi_{1,0}(s)=1/2}$ is
``ergodizd'' as time passes from $t_{k,0}$ to $t_{k+1,0}$.  If
$\varepsilon$ is sufficiently small and $t_{k+1,\varepsilon}\leq
(T\wedge T_\varepsilon)/\varepsilon$, then the spacings between
these times are uniformly of order $1$, i.e.~$1/\omega
<t_{k+1,\varepsilon}-t_{k,\varepsilon}< \omega$.  Thus,
\begin {equation}\label{eq:hc_splitting}
\begin {split}
    \sup_{0\leq \tau\leq
    T\wedge T_\varepsilon}
    &
    \abs{\int_0^{\tau/\varepsilon}
    s_{1,\varepsilon}(s)\delta_{\varphi_{1,\varepsilon}(s)=1/2}-
    \frac{s_{1,\varepsilon}(s)^2}{2X_\varepsilon(s)}ds}
    \\
    &
    \leq
    \mathcal{O} (1)+
    \sum_{t_{k+1,\varepsilon}\leq
    \frac {T\wedge T_\varepsilon}{\varepsilon}}
    \abs{
    \int_{t_{k,\varepsilon}}^{t_{k+1,\varepsilon}}
    s_{1,\varepsilon}(s)\delta_{\varphi_{1,\varepsilon}(s)=1/2}-
    \frac{s_{1,\varepsilon}(s)^2}{2X_\varepsilon(s)}ds}.
\end {split}
\end {equation}

\paragraph*{Step 4:  Control of individual terms by comparison with
        solutions along fibers.}

The sum in Equation \eqref{eq:hc_splitting} has no more than
$\mathcal{O} (1/\varepsilon)$ terms, and so it suffices to show that
each term is no larger than $\mathcal{O} (\varepsilon) $.  We can
accomplish this by comparing the motions of $z_\varepsilon(t)$ for
$t_{k,\varepsilon}\leq t\leq t_{k+1,\varepsilon}$ with the solution
of the $\varepsilon=0$ version of Equation \eqref{eq:1dode} that, at
time $t_{k,\varepsilon}$, is located at
$z_\varepsilon(t_{k,\varepsilon})$. Since each term in the sum has
the same form, without loss of generality we will only examine the
first term and suppose that $t_{0,\varepsilon} =0$, i.e.~that
$\varphi_{1,\varepsilon} (0) =0$.

\begin {lem}
\label {lem:1d_divergence}

If $t_{1,\varepsilon}\leq\frac{T\wedge T_\varepsilon}{\varepsilon}$,
then
$
    \sup_{0\leq t\leq t_{1,\varepsilon}}\abs{z_{0}(t)-
    z_\varepsilon(t)}=\mathcal{O}(\varepsilon ).
$
\end {lem}
\begin {proof}
To check that $\sup_{0\leq t\leq t_{1,\varepsilon}}\abs{h_{0}(t)-
    h_\varepsilon(t)} = \mathcal{O}(\varepsilon) $,
first note that $h_0 (t) =h_0(0) =h_\varepsilon (0) $. Then
$dX_\varepsilon/dt =\mathcal{O} (\varepsilon) $, so that $X_0(t)
-X_\varepsilon(t) =\mathcal{O} (\varepsilon t) $.  Furthermore, the
other slow variables change by $\mathcal{O} (\varepsilon) $ at
collisions, while the number of collisions in the time interval $[0,
t_{1,\varepsilon}] $ is $ \mathcal{O} (1) $.

It remains to show that $\sup_{0\leq t\leq t_{1,\varepsilon}}
\abs{\varphi_{i,0}(t)-\varphi_{i,\varepsilon}(t)}=
\mathcal{O}(\varepsilon )$. Using what we know about the divergence
of the slow variables,
\[
\begin {split}
    \varphi_{1,0} (t) -\varphi_{1,\varepsilon} (t)  &=
    \int_0^t\frac{s_{1,0}(s)}{2X_0(s)}-
    \frac{s_{1,\varepsilon}(s)}{2X_\varepsilon(s)}+
    \mathcal{O}(\varepsilon)ds
    =\int_0 ^t
    \mathcal{O} (\varepsilon)ds
    = \mathcal{O} (\varepsilon )\\
\end {split}
\]
for $0\leq t \leq t_{1,\varepsilon} $.  Showing that $\sup_{0\leq
t\leq t_{1,\varepsilon}}
\abs{\varphi_{2,0}(t)-\varphi_{2,\varepsilon}(t)}=
\mathcal{O}(\varepsilon )$ is similar.

\end {proof}

From Lemma \ref{lem:1d_divergence},
$t_{1,\varepsilon}=t_{1,0}+\mathcal{O}(\varepsilon)
=2X_0/s_{1,0}+\mathcal{O}(\varepsilon)$.  We conclude that
\[
\begin {split}
    \int_{0}^{t_{1,\varepsilon}}
    s_{1,\varepsilon}(s)\delta_{\varphi_{1,\varepsilon}(s)=1/2}-
    \frac{s_{1,\varepsilon}(s)^2}{2X_\varepsilon(s)}ds
    & =
    \mathcal{O} (\varepsilon) +
    \int_{0}^{t_{1,\varepsilon}}
    s_{1,0}(s)\delta_{\varphi_{1,\varepsilon}(s)=1/2}-
    \frac{s_{1,0}(s)^2}{2X_0(s)}ds
    \\
    &=
    \mathcal{O}(\varepsilon)+s_{1,0}
    -t_{1,\varepsilon}\frac{s_{1,0}^2}{2X_0}
    =\mathcal{O}(\varepsilon).
\end {split}
\]

From Equations \eqref{eq:1dgronwall} and \eqref{eq:hc_splitting}, we
see that $
    \sup_{0\leq \tau\leq T\wedge T_\varepsilon}
    \abs{h_\varepsilon(\tau/\varepsilon)-\bar
    h(\tau)}
    =\mathcal{O}(\varepsilon ),
$ independent of the initial condition in $h^ {-1}\mathcal{V}$.

\subsection{Extension to multiple gas particles}

When $n_1,n_2>1$, only minor modifications are necessary to
generalize the proof above. We start by extending the slow variables
$h$ to a full set of coordinates on phase space by defining the
angle variables $\varphi _{i,j}\in [0,1]/ \, 0\sim 1=S^1$ for $1\leq
i\leq 2, $ $1\leq j\leq n_i$:
\[
\begin {split}
    \varphi_{1,j} =\varphi_{1,j} (x_{1,j},v_{1,j},X)=
    &\begin {cases}
    \frac{x_{1,j}}{2X} &\text { if } v_{1,j} >0\\
    1-\frac{x_{1,j}}{2X} &\text { if } v_{1,j}<0\\
    \end {cases}
    \\
    \varphi_{2,j}=\varphi_{2,j}(x_{2,j},v_{2,j},X)=
    &\begin {cases}
    \frac{1-x_{2,j}}{2(1-X)} &\text { if } v_{2,j}<0\\
    1- \frac{1-x_{2,j}}{2(1-X)} &\text { if } v_{2,j}>0\\
    \end {cases}\\
\end {split}.
\]
Then $d\varphi_{1,j}/dt=s_{1,j}(2X) ^ {-1} +\mathcal{O}
(\varepsilon)$, $d\varphi_{2,j}/dt=s_{2,j}(2(1-X)) ^ {-1}
+\mathcal{O} (\varepsilon)$, and $z =
(h,\varphi_{1,j},\varphi_{2,j}) $ represents a choice of coordinates
on our phase space, which is homeomorphic to $(\text{a subset of }
\mathbb{R}^{n_1+n_2+2})\times \mathbb{ T} ^ {n_1+n_2} $.  In these
coordinates, the dynamical system yields a discontinuous flow
$z_\varepsilon(t) $ on phase space.  The flow preserves Liouville
measure, which in our coordinates has a density proportional to $X^
{n_1}(1-X)^{n_2}$. As is Section \ref{sct:1dps_hcPS}, one can show
that the measure of initial conditions leading to multiple particle
collisions is zero.

Next, define $H(z) $ by
\[
    H(z) =
    \begin{bmatrix}
    W\\
    \sum_{j=1}^{n_1}2m_{1,j} s_{1,j} \delta_{\varphi_{1,j}=1/2}
    -\sum_{j=1}^{n_2}2m_2
     s_{2,j}\delta_{\varphi_{2,j}=1/2}\\
    -2W\delta_{\varphi_{1,j}=1/2}\\
    2W\delta_{\varphi_{2,j}=1/2}\\
    \end{bmatrix}.
\]
For $0 \leq t\leq \frac{T\wedge T_\varepsilon}{\varepsilon}$, $
    h_\varepsilon(t)-h_\varepsilon(0)=
    \varepsilon\int_0^t
    H(z_\varepsilon(s))ds+\mathcal{O}(\varepsilon).
$
From here, the rest of the proof follows the same arguments made in
Section \ref{sct:1dps_hc_unif}.

\section{Proof of the main result for the soft core piston}
\label {sct:sc_proof}

For the remainder of this work, we consider the family of
Hamiltonian systems introduced in Section
\ref{sct:soft_core_results}, which are parameterized by
$\varepsilon,\delta\geq 0$. For simplicity, we specialize to
$n_1=n_2=1$. As in Section \ref{sct:hc_proof}, the generalization to
$n_1, n_2>1$ is not difficult. The Hamiltonian dynamics are given by
the following ordinary differential equation:
\begin {equation}
\label {eq:smooth1dode1}
\begin {split}
    \frac{dX}{dt}&=\varepsilon W,\\
    \frac{dW}{dt}&=\varepsilon\left(
    -\kappa_\delta'(X-x_{1})
    +\kappa_\delta'(x_{2}-X)\right),\\
    \frac{dx_{1}}{dt}&=v_{1},\\
    \frac{dv_{1}}{dt}&=\frac{1}{m_{1}}\bigl( -\kappa_\delta'(x_{1})
    +\kappa_\delta'(X-x_{1})\bigr),\\
    \frac{dx_{2}}{dt}&=v_{2},\\
    \frac{dv_{2}}{dt}&=\frac{1}{m_{2}}\bigl( -\kappa_\delta'(x_{2}-X)
    +\kappa_\delta'(1-x_{2})\bigr).\\
\end {split}
\end {equation}
Recalling the particle energies defined by Equation
\eqref{eq:soft_energies}, we find that
\begin {equation*}
\begin {split}
    \frac {d E_{1}} {dt} =\varepsilon W\kappa_\delta' (X
    -x_{1}),\qquad
    \frac {d E_{2}} {dt} = -\varepsilon W\kappa_\delta' (x_{2} -X).\\
\end {split}
\end {equation*}

\comment {
\begin {equation}
\label {eq:smooth1dode1}
\begin {split}
    \frac{dX}{dt}&=\varepsilon W,\\
    \frac{dW}{dt}&=\varepsilon\left(\sum_{j=1}^{n_1}
    -\kappa_\delta'(X-x_{1,j})
    +\sum_{j=1}^{n_2}\kappa_\delta'(x_{2,j}-X)\right),\\
    \frac{dx_{1,j}}{dt}&=v_{1,j},\\
    \frac{dv_{1,j}}{dt}&=\frac{1}{m_{1,j}}\bigl( -\kappa_\delta'(x_{1,j})
    +\kappa_\delta'(X-x_{1,j})\bigr),\\
    \frac{dx_{2,j}}{dt}&=v_{2,j},\\
    \frac{dv_{2,j}}{dt}&=\frac{1}{m_{2,j}}\bigl( -\kappa_\delta'(x_{2,j}-X)
    +\kappa_\delta'(1-x_{2,j})\bigr).\\
\end {split}
\end {equation}
Recalling the particle energies $E_{i,j} $ defined by Equation
\eqref{eq:soft_energies}, we find that
\begin {equation*}
\begin {split}
    \frac {d E_{1,j}} {dt} =\varepsilon W\kappa_\delta' (X
    -x_{1,j}),\qquad
    \frac {d E_{2,j}} {dt} = -\varepsilon W\kappa_\delta' (x_{2,j} -X).\\
\end {split}
\end {equation*}

}

For the compact set $\mathcal{V} $ introduced in Section
\ref{sct:1d_sm_results}, fix a small positive number $\mathcal{E} $
and an open set $\mathcal{U}\subset \mathbb{R}^4 $ such that $
\mathcal{V} \subset\mathcal{U} $ and $h\in\mathcal{U}\Rightarrow
X\in (\mathcal{E},1 -\mathcal{E})$, $W\subset\subset \mathbb {R} $,
and $\mathcal{E} <E_1,E_2 <\kappa (0) -\mathcal{E} $. We only
consider the dynamics for $0<\delta<\mathcal{E}/2$ and $h
\in\mathcal{U} $.

Define
\begin {equation*}
\begin {split}
    U_1(x_1)&=U_1(x_1,X,\delta)=\kappa_\delta(x_1)
    +\kappa_\delta(X-x_1),\\
    U_2(x_2)&=U_2(x_2,X,\delta)  =\kappa_\delta(x_2-X)
    +\kappa_\delta(1-x_2).\\
\end {split}
\end {equation*}
Then the energies $E_i$ satisfy $E_i=m_i v_i^2/2+U_i(x_i)$.

Let $T_1=T_1(X,E_1,\delta) $ and $T_2=T_2(X,E_2,\delta)$ denote the
periods of the motions of the left and right gas particles,
respectively, when $\varepsilon =0 $.

\begin {lem}
\label {lem:1d_smooth_periods}
For $i=1,2 $,
\[
    T_i\in\mathcal{C} ^1\{(X,E_i,\delta):
    X\in (\mathcal{E},1-\mathcal{E})
    ,E_i\in (\mathcal{E} ,\kappa (0) -\mathcal{E}),0\leq
    \delta <\mathcal{E}/2\} .
\]
Furthermore,
\[
\begin {split}
    T_1(X,E_1,\delta)&=\sqrt{\frac  {2m_1} {E_1}}X
    +\mathcal{O} (\delta),\quad
    T_2(X,E_2,\delta)=\sqrt{\frac  {2m_2} {E_2}} (1 - X)
    +\mathcal{O} (\delta).\\
\end {split}
\]
\end {lem}

The proof of this lemma is mostly computational, and so we delay it
until Section \ref{sct:smooth1DTech}.  Note especially that the
periods can be suitably defined such that their regularity extends
to $\delta=0$.

In this section, and in Section \ref{sct:smooth1DTech} below, we
adopt the following convention on the use of the $\mathcal{O} $
notation.  \emph{All use of the $\mathcal{O} $ notation will
explicitly contain the dependence on $\varepsilon$ and $\delta$} as
$\varepsilon, \delta\rightarrow 0$.  For example, if a function
$f(h,\varepsilon,\delta) =\mathcal{O} (\varepsilon) $, then there
exists $\delta',\varepsilon' >0$ such that $\sup_{0<\varepsilon\leq
\varepsilon',\,0<\delta\leq\delta',\,h\in\mathcal{V}}
\abs{f(h,\varepsilon,\delta)/\varepsilon} < \infty $.

When $\varepsilon =0 $, $\frac { dx_i} {dt} =\pm \sqrt {\frac {2}
{m_i} (E_i -U_i (x_i))} $.  Define $a = a (E_i,\delta) $ by
$\kappa_\delta (a) =\kappa (a/\delta) =E_i $, so that $a(E_1,\delta)
$ is a turning point for the left gas particle. Then $a=\delta\kappa
^ {-1} (E_i) $, where $\kappa ^ {-1} $ is defined as follows:
$\kappa: [0,1]\rightarrow [0,\kappa (0)] $ takes $0 $ to $\kappa (0)
$ and $1 $ to $0 $. Furthermore, $\kappa\in\mathcal{C} ^2 ([0,1]) $,
$\kappa '\leq 0 $, and $\kappa ' (x) <0 $ if $x < 1 $. By
monotonicity, $\kappa ^ {-1}\colon [0,\kappa (0)]\rightarrow [0,1] $
exists and takes $0 $ to $1 $ and $\kappa (0) $ to $0 $. Also, by
the Implicit Function Theorem, $\kappa ^ {-1}\in\mathcal{C} ^2
((0,\kappa (0)]) $, $(\kappa ^ {-1})'(y) < 0 $ for $y>0$, and
$(\kappa ^ {-1})'(y)\rightarrow -\infty $ as $y\rightarrow 0^ +$.
Because we only consider energies $E_i\in
(\mathcal{E},\kappa(0)-\mathcal{E}) $, it follows that $a
(E_i,\delta)$ is a $\mathcal{C} ^2$ function for the domains of
interest.

\subsection{Derivation of the averaged equation}

\label {sct:smooth_1D_derivation}

As we previously pointed out, for each fixed $\delta$, Anosov's
Theorem applies directly to the family of ordinary differential
equations in Equation \eqref{eq:smooth1dode1}, provided that
$\delta$ is sufficiently small. The invariant fibers $\mathcal{M}_h$
of the $\varepsilon=0$ flow are tori described by a fixed value of
the four slow variables and $\{(X,W,x_1,v_1,x_2,v_2): E_1=m_1
v_1^2/2+U_1(x_1,X,\delta), E_2=m_2 v_2^2/2+U_2(x_2,X,\delta)\} $. If
we use $(x_1,x_2) $ as local coordinates on $\mathcal{M}_h$, which
is valid except when $v_1\text { or }v_2=0$, the invariant measure
$\mu_h$ of the unperturbed flow has the density
\[
    \frac {dx_1 dx_2}
    {
    T_1\sqrt {\frac {2} {m_1} (E_1-U_1(x_1))}\:
    T_2\sqrt {\frac {2} {m_2} (E_2-U_2(x_2))}
    }.
\]
The restricted flow is ergodic for almost every $h $. See Corollary
\ref{cor:irrat_periods} in Section \ref{sct:smooth1DTech}.

Now
\[
    \frac {dh_\varepsilon ^\delta} {dt}
    =\varepsilon
    \begin {bmatrix}
             W \\
         -\kappa_\delta'(X -x_1 )
                +\kappa_\delta'(x_2-X )\\
             W\kappa_\delta' (X
                -x_1)\\
            - W\kappa_\delta' (
                x_2-X )\\
    \end {bmatrix},
\]
and
\[
\begin {split}
    \int_{\mathcal{M}_h}\kappa_\delta'(X -x_1 )d\mu_h
    &
    =\frac {2} {T_1}\int_a^{X-a}dx_1\frac{\kappa_\delta'(X-x_1)}
    {\sqrt {\frac {2} {m_1} (E_1-U_1(x_1))}}
    =\frac {\sqrt{2m_1}} {T_1}\int_{X-\delta}^{X-a}dx_1\frac{\kappa_\delta'(X-x_1)}
    {\sqrt { E_1-\kappa_\delta(X-x_1)}}
    \\
    &
    =-\frac {\sqrt{2m_1}} {T_1}\int_{0}^{E_1}\frac{du}
    {\sqrt { E_1-u}}
    = -\frac {\sqrt{8m_1E_1}} {T_1}.
    \\
\end {split}
\]
Similarly,
\[
\begin {split}
    \int_{\mathcal{M}_h}\kappa_\delta'(x_2-X)d\mu_h
    = -\frac {\sqrt{8m_2E_2}} {T_2}.\\
\end {split}
\]
It follows that the averaged vector field is
\[
    \bar H^\delta (h) =
    \begin {bmatrix}
             W \\
         \frac{\sqrt{8m_1E_1}}{T_1}-
     \frac{\sqrt{8m_2E_2}}{T_2}\\
             -W\frac{\sqrt{8m_1E_1}}{T_1}\\
            +W\frac{\sqrt{8m_2E_2}}{T_2}\\
        \end {bmatrix},
\]
where from Lemma \ref{lem:1d_smooth_periods} we see that $\bar H^
{\cdot} (\cdot)\in\mathcal{C} ^1 (\{(\delta,h):0\leq\delta
<\mathcal{E}/2,h\in\mathcal{V}\})$.  $\bar H^0(h) $ agrees with the
averaged vector field for the hard core system from Equation
\eqref{eq:1davg}, once we account for the change of coordinates
$E_i=m_i s_i^2/2$.

\begin {rem}

An argument due to Neishtadt and Sinai~\cite{NS04} shows that the
solutions to the averaged equation \eqref{eq:smooth_1D_averaged_eq}
are periodic.  This argument also shows that, as in the case $\delta
=0$, the limiting dynamics of $(X,W) $ are effectively Hamiltonian,
with the shape of the Hamiltonian depending on $\delta$, $X(0) $,
and the initial energies of the gas particles. The argument depends
heavily on the observation that the phase integrals
\[
    I_i(X, E_i,\delta)=
    \int_{\frac{1}{2}m_i v^2+U_i(x,X,\delta)\leq E_i}
    dxdv
\]
are adiabatic invariants, i.e.~they are integrals of the solutions
to the averaged equation.  Thus the four-dimensional phase space of
the averaged equation is foliated by invariant two-dimensional
submanifolds, and one can think of the effective Hamiltonians for
the piston as living on these submanifolds.

\end {rem}

\subsection {Proof of Theorem \ref{thm:1D_smooth_uniform}}

The following arguments are motivated by our proof in Section
\ref{sct:hc_proof}, although the details are more involved as we
show that the rate of convergence is independent of all small
$\delta$.

\subsubsection{A choice of coordinates on phase space}

We wish to describe the dynamics in a coordinate system inspired by
the one used in Section \ref{sct:1dps_hcPS}.  For each fixed
$\delta\in (0,\delta_0]$, this change of coordinates will be
$\mathcal{C} ^1 $ in all variables on the domain of interest.
However, it is an exercise in analysis to show this, and so we delay
the proofs of the following two lemmas until Section
\ref{sct:smooth1DTech}.

We introduce the angular coordinates $\varphi _i\in [0,1]/ \, 0\sim
1=S^1$ defined by
\begin {equation}
\label {eq:anglevar_defn}
\begin {split}
    \varphi_1 =\varphi_1 (x_1,v_1,X)=&\begin {cases}
    0&\text { if } x_1 =a\\
    \frac{1}{T_1} \int_a^{x_1}\sqrt{\frac{m_1/2}
        {E_1-U_1(s)}}ds&\text { if } v_1 >0\\
    1/2&\text { if } x_1 =X -a\\
    1-\frac{1}{T_1} \int_a^{x_1}\sqrt{\frac{m_1/2}
        {E_1-U_1(s)}}ds &\text { if } v_1<0\\
    \end {cases}
    \\
    \varphi_2=\varphi_2(x_2,v_2,X)=&\begin {cases}
    0&\text { if } x_2 =1-a\\
    \frac{1}{T_2} \int_{x_2} ^{1 - a}\sqrt{\frac{m_2/2}
        {E_2-U_2(s)}}ds&\text { if } v_2 < 0\\
    1/2&\text { if } x_2 =X+a\\
    1-\frac{1}{T_2} \int_{x_2} ^{1 -a}\sqrt{\frac{m_2/2}
        {E_2-U_2(s)}}ds &\text { if } v_2 > 0\\
    \end {cases}\\
\end {split}.
\end{equation}
Then $z =(h,\varphi_1,\varphi_2) $ is a choice of coordinates on $h^
{-1}\mathcal{U} $. As before, we will abuse notation and let $h (z)
$ denote the projection onto the first four coordinates of $z $.

There is a fixed value of $\delta_0$ in the statement of Theorem
\ref{thm:1D_smooth_uniform}.  However, for the purposes of our
proof, it will be convenient to progressively choose $\delta_0$
smaller when needed.  At the end of the proof, we will have only
shrunk $\delta_0$ a finite number of times, and this final value
will satisfies the requirements of the theorem.  Our first
requirement on $\delta_0$ is that it is smaller than
$\mathcal{E}/2$.

\begin {lem}
\label {lem:1d_smooth_ode}

If $\delta_0 >0$ is sufficiently small, then for each $\delta\in
(0,\delta_0]$ the ordinary differential equation
\eqref{eq:smooth1dode1} in the coordinates $z $ takes the form
\begin {equation}
\label {eq:smooth1dode_2}
    \frac {dz } {dt} =Z ^\delta (z ,\varepsilon),
\end {equation}
where $Z ^\delta \in\mathcal{C} ^1 ( h  ^ {-1} \mathcal{U}\times
[0,\infty)) $. When $z \in h^{-1}\mathcal{U} $,

\begin {equation}
\label {eq:smooth1dode_3}
    Z ^\delta (z,\varepsilon) =
    \begin {bmatrix}
            \varepsilon W \\
            \varepsilon\bigl
                ( -\kappa_\delta'(X -x_1 (z))
                +\kappa_\delta'(x_2(z)-X )\bigr)\\
            \varepsilon W\kappa_\delta' (X
                -x_1(z))\\
            -\varepsilon W\kappa_\delta' (
                x_2(z)-X )\\
            \frac {1} {T_1}+\mathcal{O}
                (\varepsilon)\\
            \frac {1} {T_2}+\mathcal{O}
                (\varepsilon)\\
    \end {bmatrix}.
\end {equation}

\end {lem}

Recall that, by our conventions, the $\mathcal{O} (\varepsilon) $
terms in Equation \eqref{eq:smooth1dode_3} have a size that can be
bounded independent of all $\delta$ sufficiently small.  Denote the
flow determined by $Z ^\delta (\cdot,\varepsilon) $ by
$z_\varepsilon ^\delta (t) $, and its components by
$X_\varepsilon^\delta (t) $, $W_\varepsilon^\delta (t) $,
$E_{1,\varepsilon}^\delta (t) $, etc. Also, set
$h_\varepsilon^\delta (t) =h(z_\varepsilon ^\delta (t) )$.  From
Equation \eqref{eq:smooth1dode_3},
\begin{equation}\label {eq:sm_H}
    H ^\delta (z,\varepsilon) =
    \frac{1}{\varepsilon}
    \frac{dh_\varepsilon^\delta}{dt}=
    \begin {bmatrix}
            W \\
             -\kappa_\delta'(X -x_1 (z))
                +\kappa_\delta'(x_2(z)-X )\\
            W\kappa_\delta' (X
                -x_1(z))\\
            -W\kappa_\delta' (
                x_2(z)-X )\\
    \end {bmatrix}.
\end {equation}
In particular, $H ^\delta (z,\varepsilon) =H ^\delta (z,0)$.

Before proceeding, we need one final technical lemma.

\begin {lem}
\label {lem:1d_smooth_delta}

If $\delta_0>0$ is chosen sufficiently small, there exists a
constant $K$ such that for all $\delta\in (0,\delta_0]$,
$\kappa_\delta' (\abs{X -x_i(z)}) =0 $ unless $\varphi_i\in [1/2
-K\delta,1/2 +K\delta] $.
\end {lem}

\subsubsection{Argument for uniform convergence}

We start by proving the following lemma, which essentially says that
an orbit $z_\varepsilon ^\delta (t) $ only spends a fraction
$\mathcal{O} (\delta) $ of its time in a region of phase space where
$\abs{H^\delta (z_\varepsilon ^\delta (t),\varepsilon)}
=\abs{H^\delta (z_\varepsilon ^\delta (t),0)}$ is of size
$\mathcal{O} (\delta ^ {-1}) $

\begin {lem}
\label {lem:1d_smooth_intbound}

For $0\leq\mathcal{T}'\leq\mathcal{T}\leq\frac {T\wedge
T_\varepsilon ^\delta} {\varepsilon}$,
\[
    \int_{\mathcal{T}'}^{\mathcal{T}} H ^\delta (z_\varepsilon ^\delta
    (s),0)ds
    =\mathcal{O} (1\vee(\mathcal{T}-\mathcal{T}')).
\]
\end {lem}

\begin {proof}

Without loss of generality, $\mathcal{T}' =0$.  From Lemmas
\ref{lem:1d_smooth_periods} and \ref{lem:1d_smooth_ode} it follows
that if we choose $\delta_0$ sufficiently small, then there exists
$\omega
>0 $ such that for all sufficiently small $\varepsilon $ and all
$\delta\in (0,\delta_0]$, $h\in\mathcal{V}\Rightarrow 1/\omega
<\frac {d\varphi_{i,\varepsilon} ^\delta} {dt} <\omega $.   Define
the set $B = [1/2 -K\delta, 1/2+K\delta] $, where $K $ comes from
Lemma \ref{lem:1d_smooth_delta}.  Then we find a crude bound on
$\int_0^{\mathcal{T}} \abs {\kappa_\delta' \bigl(X_\varepsilon
^\delta (s) -x_1 (z_\varepsilon ^ \delta(s))\bigr)}ds $ using that
\[
    \frac {d\varphi_{1,\varepsilon} ^\delta} {dt}\text { is }
    \begin {cases}
    \geq 1/\omega &\text { if } \varphi_{1,\varepsilon} ^\delta\in B\\
    \leq \omega &\text { if } \varphi_{1,\varepsilon} ^\delta\in B^c.\\
    \end {cases}
\]
This yields
\[
\begin{split}
    \int_0^{\mathcal{T}}  \abs {\kappa_\delta'
    \bigl( X_\varepsilon ^\delta(s)
    -x_1
    (z_\varepsilon ^ \delta(s))  \bigr)}ds
    &
    \leq\frac { \C} {\delta}\int_0^{\mathcal{T}}
    1_{\varphi_1 (s)\in B}ds\\
    & \leq \frac { \C} {\delta} \left(\frac {2K\omega\delta}
    {2K\omega\delta+\frac{1-2K\delta}{\omega}}\mathcal{T}
    +2K\omega\delta\right)=\mathcal{O}(1\vee\mathcal{T}).\\
\end{split}
\]
Similarly, $\int_0^{\mathcal{T}} \abs {\kappa_\delta' (
x_2(z_\varepsilon ^ \delta(s))-X_\varepsilon ^\delta(s))}ds
    =\mathcal{O}(1\vee\mathcal{T})$, and so $\int_0^{\mathcal{T}}
    H ^\delta (z_\varepsilon ^\delta
    (s),0)ds= \mathcal{O}(1\vee\mathcal{T})$.
\end {proof}

We now follow steps one through four from Section
\ref{sct:1dps_hc_unif}, making modifications where necessary.

\paragraph*{Step 1:  Reduction using Gronwall's Inequality.}

Now $h_\varepsilon^\delta(\tau/\varepsilon)$ satisfies $
h_\varepsilon^\delta(\tau/\varepsilon)-h_\varepsilon^\delta(0)    =
\varepsilon\int_0^{\tau/\varepsilon}
H^\delta(z_\varepsilon^\delta(s),0)ds $.  Define
$e_\varepsilon^\delta(\tau) =\varepsilon\int_0^{\tau/\varepsilon}
H^\delta(z_\varepsilon^\delta(s),0)- \bar
H^\delta(h_\varepsilon^\delta(s))ds$. It follows from Gronwall's
Inequality and the fact that $\bar H^ {\cdot} (\cdot)\in\mathcal{C}
^1 (\{(\delta,h):0\leq\delta \leq\delta_0,h\in\mathcal{V}\})$ that
\begin {equation}
\label{eq:1d_smooth_gronwall}
    \sup_{0\leq \tau\leq T\wedge T_\varepsilon ^\delta}
    \abs{h_\varepsilon ^\delta (\tau/\varepsilon)-\bar
    h ^\delta (\tau)} \leq
    \left(\sup_{0\leq \tau\leq T\wedge T_\varepsilon ^\delta}
    \abs{e_\varepsilon^\delta (\tau)}\right)
    e^{ \Lip{\bar
    H ^\delta \arrowvert _\mathcal{V}} T}
    =\mathcal{O}\left(\sup_{0\leq \tau\leq T\wedge T_\varepsilon ^\delta}
    \abs{e_\varepsilon^\delta (\tau)}\right).
\end {equation}

\paragraph*{Step 2:  A splitting according to particles.}

Next,
\[
    H^\delta(z,0)-\bar H^\delta(h) =
    \begin{bmatrix}
    0\\
    -\kappa_\delta'(X-x_1(z))-\frac {\sqrt{8m_1 E_1}} {T_1}\\
    W\kappa_\delta'(X-x_1(z))+W\frac {\sqrt{8m_1 E_1}} {T_1}\\
    0\\
    \end{bmatrix}
    +
    \begin{bmatrix}
    0\\
    \kappa_\delta'(x_2(z)-X)+\frac {\sqrt{8m_2 E_2}} {T_2}\\
    0\\
    -W\kappa_\delta'(x_2(z)-X)-W\frac {\sqrt{8m_2 E_2}} {T_2}\\
    \end{bmatrix}
    ,
\]
and so, in order to show that $\sup_{0\leq \tau\leq T\wedge
T_\varepsilon^\delta}\abs{e_\varepsilon^\delta(\tau)} =\mathcal{O}
(\varepsilon) $, it suffices to show that for $i=1,2$,
\begin {equation*}
\begin {split}
    \sup_{0\leq \tau\leq
    T\wedge T_\varepsilon^\delta}
    &
    \abs{\int_0^{\tau/\varepsilon}
    \kappa_\delta'\bigl(\abs{X_\varepsilon^\delta(s)-
    x_i(z_\varepsilon^\delta(s))}\bigr)+
    \frac {\sqrt{8m_i E_{i,\varepsilon} ^\delta(s)}} {T_i(X_{\varepsilon}^\delta(s),
     E_{i,\varepsilon}^\delta(s),
    \delta)}ds}
    =
    \mathcal{O} (1),
    \\
    \sup_{0\leq \tau\leq
    T\wedge T_\varepsilon^\delta}
    &
    \abs{\int_0^{\tau/\varepsilon}W_\varepsilon(s)
    \kappa_\delta'\bigl(\abs{X_\varepsilon^\delta(s)-
    x_i(z_\varepsilon^\delta(s))}\bigr)+W_\varepsilon(s)
    \frac {\sqrt{8m_i E_{i,\varepsilon} ^\delta(s)}} {T_i(X_{\varepsilon}^\delta(s),
     E_{i,\varepsilon}^\delta(s),
    \delta)}ds}
    =
    \mathcal{O} (1).\\
\end {split}
\end {equation*}
We only demonstrate that $    \sup_{0\leq \tau\leq
    T\wedge T_\varepsilon^\delta}
    \abs{\int_0^{\tau/\varepsilon} \kappa_\delta'\bigl(X_\varepsilon^\delta(s)-
    x_1(z_\varepsilon^\delta(s))\bigr)+
    \frac {\sqrt{8m_1 E_{1,\varepsilon} ^\delta(s)}} {T_1(X_{\varepsilon}^\delta(s),
     E_{1,\varepsilon}^\delta(s),
    \delta)}ds}
    =
    \mathcal{O} (1)$.
The other three terms are handled similarly.

\paragraph*{Step 3:  A sequence of times adapted for ergodization.}

Define the sequence of times $t_{k,\varepsilon}^\delta $ inductively
by $t_{0,\varepsilon}^\delta=\inf\{t\geq
0:\varphi_{1,\varepsilon}^\delta(t)=0\}$,
$t_{k+1,\varepsilon}^\delta=\inf\{t>t_{k,\varepsilon}^\delta:
\varphi_{1,\varepsilon}^\delta(t)=0\}$. If $\varepsilon$ and
$\delta$ are sufficiently small and $t_{k+1,\varepsilon}^\delta\leq
(T\wedge T_\varepsilon^\delta)/\varepsilon$, then it follows from
Lemma \ref{lem:1d_smooth_ode} that $1/\omega
<t_{k+1,\varepsilon}^\delta-t_{k,\varepsilon}^\delta< \omega$. From
Lemmas \ref{lem:1d_smooth_ode} and \ref{lem:1d_smooth_intbound} it
follows that
\begin {equation}\label{eq:sc_splitting}
\begin {split}
    &\sup_{0\leq \tau\leq
    T\wedge T_\varepsilon^\delta}
    \abs{\int_0^{\tau/\varepsilon} \kappa_\delta'\bigl(X_\varepsilon^\delta(s)-
    x_1(z_\varepsilon^\delta(s))\bigr)+
    \frac {\sqrt{8m_1 E_{1,\varepsilon} ^\delta(s)}} {T_1(X_{\varepsilon}^\delta(s),
     E_{1,\varepsilon}^\delta(s),
    \delta)}ds}
    \\
    &\leq
    \mathcal{O} (1)+
    \sum_{t_{k+1,\varepsilon}^\delta\leq
    \frac {T\wedge T_\varepsilon^\delta}{\varepsilon}}
    \abs{\int_{t_{k,\varepsilon}^\delta}^{t_{k+1,\varepsilon}^\delta}
     \kappa_\delta'\bigl(X_\varepsilon^\delta(s)-
    x_1(z_\varepsilon^\delta(s))\bigr)+
    \frac {\sqrt{8m_1 E_{1,\varepsilon} ^\delta(s)}} {T_1(X_{\varepsilon}^\delta(s),
    E_{1,\varepsilon}^\delta(s),
    \delta)}ds}.
\end {split}
\end {equation}

\paragraph*{Step 4:  Control of individual terms by comparison with
        solutions along fibers.}

As before, it suffices to show that each term in the sum in Equation
\eqref{eq:sc_splitting} is no larger than $\mathcal{O} (\varepsilon)
$. Without loss of generality we will only examine the first term
and suppose that $t_{0,\varepsilon}^\delta =0$, i.e.~that
$\varphi_{1,\varepsilon}^\delta (0) =0$.

\begin {lem}
\label {lem:1d_smooth_divergence}

If $t_{1,\varepsilon}^\delta\leq\frac{T\wedge
T_\varepsilon^\delta}{\varepsilon}$, then
$
    \sup_{0\leq t\leq t_{1,\varepsilon}^\delta}\abs{z_{0}^\delta(t)-
    z_\varepsilon^\delta(t)}=\mathcal{O}(\varepsilon ).
$
\end {lem}
\begin {proof}
By Lemma \ref{lem:1d_smooth_intbound}, $h_{0} ^\delta
(t)-h_\varepsilon ^\delta (t) =h_{\varepsilon} ^\delta
(0)-h_\varepsilon ^\delta (t) = -\varepsilon\int_0 ^t H^\delta
(z_{\varepsilon} ^\delta (s),0)ds=\mathcal{O} (\varepsilon (1\vee
t))$ for $t\geq 0$.

Using what we know about the divergence of the slow variables, we
find that
\[
\begin {split}
    \varphi_{1,0} ^\delta (t) -\varphi_{1,\varepsilon} ^\delta (t)
    &=
    \int_0^t\frac{1}{T_1 (X_0 ^\delta (s),E_0 ^\delta (s),\delta)}-
    \frac{1}{T_1 (X_\varepsilon ^\delta (s),E_\varepsilon ^\delta
    (s),\delta)}+
    \mathcal{O}(\varepsilon)ds\\
    &=\int_0 ^t
    \mathcal{O} (\varepsilon)ds
    = \mathcal{O} (\varepsilon )\\
\end {split}
\]
for $0\leq t \leq t_{1,\varepsilon}^\delta$. Lemmas
\ref{lem:1d_smooth_periods} and \ref{lem:1d_smooth_ode} ensure the
desired uniformity in the sizes of the orders of magnitudes.
 Showing that $\sup_{0\leq t\leq t_{1,\varepsilon}^\delta} \abs{\varphi_{2,0}
^\delta (t)-\varphi_{2,\varepsilon} ^\delta (t)}=
\mathcal{O}(\varepsilon )$ is similar.

\end {proof}

From Lemma \ref{lem:1d_smooth_divergence} we find that
$t_{1,\varepsilon}=t_{1,0}+\mathcal{O}(\varepsilon) =T_1 (X_0
^\delta,E_0 ^\delta,\delta)+\mathcal{O}(\varepsilon)$.  Hence
\[
\begin {split}
    \int_{0}^{t_{1,\varepsilon}^\delta}
    \frac {\sqrt{8m_1 E_{1,\varepsilon} ^\delta(s)}} {T_1(X_{\varepsilon}^\delta(s),
    E_{1,\varepsilon}^\delta(s),
    \delta)}ds
    &=
    \mathcal{O} (\varepsilon) +
    \int_{0}^{t_{1,0}^\delta}
    \frac {\sqrt{8m_1 E_{1,0} ^\delta}} {T_1(X_{0}^\delta,
    E_{1,0}^\delta,
    \delta)}ds
     =\mathcal{O} (\varepsilon) +\sqrt{8m_1 E_{1,0} ^\delta}.
\end {split}
\]
But when $x_1(z_\varepsilon^\delta)<X_\varepsilon^\delta-a$,
\[
\begin {split}
    \frac {d} {ds}
    &
    \sqrt{E_{1,\varepsilon}^\delta(s)
    -\kappa_\delta\bigl(X_\varepsilon^\delta(s)
    -x_1(z_\varepsilon^\delta(s))\bigr)}
    =
    \frac {\text {sign} \bigl(v_1(z_\varepsilon^\delta(s))\bigr)
    \kappa_\delta ' \bigl(X_\varepsilon^\delta(s)
        -x_1(z_\varepsilon^\delta(s))\bigr)} {\sqrt{2m_1}},
    \\
\end {split}
\]
and so
\[
    \int_{0}^{t_{1,\varepsilon}^\delta}\kappa_\delta'\bigl(X_\varepsilon^\delta(s)-
    x_1(z_\varepsilon^\delta(s))\bigr) ds
    =-\sqrt{2m_1 E_{1,\varepsilon} ^\delta(0)}
    -\sqrt{2m_1 E_{1,\varepsilon} ^\delta(t_{1,\varepsilon}^\delta)}
    =\mathcal{O} (\varepsilon) -\sqrt{8m_1 E_{1,0} ^\delta}.
\]
Hence,
\[
    \int_{0}^{t_{1,\varepsilon}^\delta}\kappa_\delta'\bigl(X_\varepsilon^\delta(s)-
    x_1(z_\varepsilon^\delta(s))\bigr)+
    \frac {\sqrt{8m_1 E_{1,\varepsilon} ^\delta(s)}} {T_1(X_{\varepsilon}^\delta(s),
    E_{1,\varepsilon}^\delta(s),
    \delta)}ds=\mathcal{O} (\varepsilon),
\]
as desired.


\section{Appendix to Section \ref{sct:sc_proof}}
\label {sct:smooth1DTech}

\paragraph* {Proof of Lemma \ref{lem:1d_smooth_periods}:}
\begin{proof}

For $0<\delta <\mathcal{E}/2 $,
\[
\begin {split}
    T_1=T_1(X,E_1,\delta)=2\int_a^{X-a}\sqrt{\frac{m_1/2}
        {E_1-U_1(s)}}ds,\quad
    T_2=T_2(X,E_2,\delta)=2\int_{X+a}^{1-a}\sqrt{\frac{m_2/2}
        {E_2-U_2(s)}}ds.\\
\end {split}
\]
We only consider the claims about $T_1 $, and for convenience we
take $m_1 =2 $.  Then
\[
\begin {split}
    T_1(X,E_1,\delta)&=
    2\int_a^{X-a}\frac{ds}{\sqrt{E_1-U_1(s)}}=
    4\int_a^{X/2}\frac{ds}{\sqrt{E_1-\kappa_\delta(s)}}\\
    &=
    4\left(\frac{X/2-\delta}{\sqrt{E_1}}+
    \int_a^{\delta}\frac{ds}{\sqrt{E_1-\kappa_\delta(s)}}\right)
    =
    \frac{2X-4\delta}{\sqrt{E_1}}+
    4\delta\int_{\kappa^{-1}(E_1)}^1
    \frac{ds}{\sqrt{E_1-\kappa(s)}}.\\
\end {split}
\]

Define
\[
    F (E) =\int_{\kappa^{-1}(E)}^1
    \frac{ds}{\sqrt{E-\kappa(s)}}
    =
    \int_{0}^E
    \frac{- (\kappa ^ {-1})' (u)}{\sqrt{E-u}}du.
\]
Notice that $ (\kappa ^ {-1})' (u)$ diverges as $u\rightarrow 0 ^ +
$, while $(E-u) ^ {-1/2} $ diverges as $u\rightarrow E ^ - $, but
both functions are still integrable on $[0,E] $.  It follows that $F
(E) $ is well defined. Then it suffices to show that $F:
[\mathcal{E},\kappa (0) -\mathcal{E}]\rightarrow \mathbb{R} $ is
$\mathcal{C} ^1 $.

Write
\[
\begin {split}
    F (E)
    = &
    \int_0 ^ {\mathcal{E}/2}
    \frac{- (\kappa ^ {-1})' (u)}{\sqrt{E-u}}du +
    \int_{\mathcal{E}/2} ^E
    \frac{- (\kappa ^ {-1})' (u)}{\sqrt{E-u}}du
    =
    F_1 (E) +F_2 (E).\\
\end {split}
\]
A standard application of the Dominated Convergence Theorem allows
us to differentiate inside the integral and conclude that
$F_1\in\mathcal{C} ^ {\infty} ([\mathcal{E} ,\kappa (0)
-\mathcal{E}]) $, with
\[
    F_1' (E) =\int_0 ^ {\mathcal{E}/2}
    \frac{ (\kappa ^ {-1})' (u)}{2 (E-u) ^ {3/2}}du.
\]

To examine $F_2 $, we make the substitution $v =E-u $ to find that
\[
    F_2 (E) =\int_0 ^ {E-\mathcal{E}/2}
    \frac{- (\kappa ^ {-1})' (E-v)}{\sqrt{v}}dv.
\]
Using the fact that $(\kappa ^ {-1})'\in\mathcal{C} ^1 ([\mathcal{E}
/2,\kappa (0)]) $ and the Dominated Convergence Theorem, we find
that $F_2 $ is differentiable, with
\[
    F_2' (E) =\frac{- (\kappa ^ {-1})' (\mathcal{E}/2)}
    {\sqrt{E-\mathcal{E}/2}} +
    \int_0 ^ {E-\mathcal{E}/2}
    \frac{- (\kappa ^ {-1})'' (E-v)}{\sqrt{v}}dv.
\]
Another application of the Dominated Convergence Theorem shows that
$F_2' $ is continuous, and so $F_2\in\mathcal{C} ^ {1} ([\mathcal{E}
,\kappa (0) -\mathcal{E}]) $.

Thus
\[
    T_1(X,E_1,\delta)=
    \frac{2X}{\sqrt{E_1}}+
    4\delta\left[-E_1 ^ {-1/2} +F_1 (E_1) +F_2 (E_1)\right]
\]
has the desired regularity.  For future reference, we note that
\begin {equation} \label {eq:partial_periods}
\begin {split}
    \frac {\partial T_1} {\partial X} =
    \frac {2} {\sqrt E_1},
    \quad
    \frac {\partial T_1} {\partial E_1} =
    \frac {-X} {E_1 ^ {3/2}} +\mathcal{O} (\delta).
\end {split}
\end {equation}

\end {proof}

\begin {cor}
\label {cor:irrat_periods} For all $\delta$  sufficiently small, the
flow $z_0 ^\delta (t) $ restricted to the invariant tori
$\mathcal{M}_c =\{h =c\} $ is ergodic (with respect to the invariant
Lebesgue measure) for almost every $c\in \mathcal{U}$.

\end {cor}

\begin {proof}
The flow is ergodic whenever the periods $T_1 $ and $T_2 $ are
irrationally related.  Fix $\delta $ sufficiently small such that $
\frac {\partial T_1} {\partial E_1} =
    -X/E_1 ^ {3/2} +\mathcal{O} (\delta) < 0 $.  Next, consider
$X $, $W$, and $E_2 $ fixed, so that $T_2 $ is constant.  Because
$T_1\in\mathcal{C}^1 $, it follows that, as we let $E_1 $ vary,
$\frac {T_1} {T_2}\notin\mathbb {Q} $ for almost every $E_1 $.  The
result follows from Fubini's Theorem.
\end {proof}

\paragraph* {Proof of Lemma \ref{lem:1d_smooth_ode}:}
\begin{proof}
For the duration of this proof, we consider the dynamics for a
small, fixed value of $\delta >0 $, which we generally suppress in
our notation.  For convenience, we take $m_1 =2 $.

Let $\psi $ denote the map taking $(X,W,x_1,v_1,x_2,v_2) $ to
$(X,W,E_1,E_2,\varphi_1,\varphi_2) $.  We claim that $\psi$ is a
$\mathcal{C} ^1 $ change of coordinates on the domain of interest.
Since $E_1 =v_1 ^2+\kappa_\delta(x_1)+\kappa_\delta(X-x_1)$, $E_1$
is a $\mathcal{C} ^2 $ function of $x_1, v_1, $ and $X$.  A similar
statement holds for $E_2 $.

The angular coordinates $\varphi_i (x_i,v_i,X) $ are defined by
Equation \eqref{eq:anglevar_defn}.  We only consider $\varphi_1 $,
as the statements for $\varphi_2 $ are similar.  Then $\varphi_1
(x_1,v_1,X) $ is clearly $\mathcal{C} ^1 $ whenever $x_1\neq a,X -a
$.  The apparent difficulties in regularity at the turning points
are only a result of how the definition of $\varphi_1 $ is presented
in Equation \eqref{eq:anglevar_defn}.  Recall that the angle
variables are actually defined by integrating the elapsed time along
orbits, and our previous definition expressed $\varphi_1 $ in a
manner which emphasized the dependence on $x_1 $.  In fact, whenever
$\abs{v_1} <\sqrt E_1 $,
\begin {equation}
\label {eq:anglevar_defn2}
\begin {split}
    \varphi_1 (x_1,v_1,X)=\begin {cases}
    -\frac{2}{T_1} \int_0^{v_1}(\kappa_\delta^{-1})'
        (E_1-v^2)dv&\text { if } x_1 <\delta\\
    \frac{1}{2}+\frac{2}{T_1} \int_0^{v_1}(\kappa_\delta^{-1})'
        (E_1-v^2)dv&\text { if } x_1 >X -\delta.\\
    \end {cases}
    \\
\end {split}
\end {equation}
Here $E_1$ is implicitly considered to be a function of $x_1, v_1,$
and $X$.  One can verify that $D\psi $ is non-degenerate on the
domain of interest, and so $\psi$ is indeed a $\mathcal{C} ^ 1 $
change of coordinates.

Next observe that $d\varphi_{1,0}/dt=1/T_1$, so Hadamard's Lemma
implies that $d\varphi_{1,\varepsilon}/dt=1/T_1 +\mathcal{O}
(\varepsilon f(\delta)) $.  It remains to show that, in fact, we may
take $f(\delta) =1$. It is easy to verify this whenever $x_1\leq X
-\delta $ because $dE_1/dt= 0 $ there.  We only perform the more
difficult verification when $x_1 > X -\delta $.

When $x_1 > X -\delta $, $\abs {v_1} <\sqrt E_1 $ and $E_1 =v_1 ^2
+\kappa_\delta (X-x_1) $.  From Equation \eqref{eq:anglevar_defn2}
we find that
\begin {equation}\label {eq:anglevar_defn3}
    \varphi_1 =    \frac{1}{2}+\frac{2\delta}{T_1(X,E_1,\delta)}
    \int_0^{v_1}(\kappa^{-1})'
        (E_1-v^2)dv.
\end {equation}
To find $d\varphi_1/dt $, we consider $\varphi_1 $ as a function of
$v_1, X, $ and $E_1 $, so that
\[
    \frac {d\varphi_1}{dt}=
    \frac {\partial \varphi_1} {\partial v_1}\frac {d v_1} {dt} +
    \frac {\partial \varphi_1} {\partial X}\frac {d X} {dt} +
    \frac {\partial \varphi_1} {\partial E_1}\frac {d E_1} {dt}.
\]
Then, using Equations \eqref{eq:partial_periods} and
\eqref{eq:anglevar_defn3}, we compute
\[
\begin {split}
    \frac {\partial \varphi_1} {\partial v_1}
    \frac {dv_1} {dt}
    & =
    \frac {2} {T_1}(\kappa_\delta^{-1})'(E_1-v_1^2)
    \frac {\kappa_\delta' (X -x_1)} {2}
    =
    \frac {1} {T_1},
    \\
    \frac {\partial \varphi_1} {\partial X}
    \frac {dX} {dt}
    & =
    \frac {1/2 -\varphi_1} {T_1}
    \frac {\partial T_1} {\partial X}
    (\varepsilon W)
    =
    \varepsilon W\frac {1/2 -\varphi_1} {T_1}
    \frac {2} {\sqrt E_1},
    \\
    \frac {\partial \varphi_1} {\partial E_1}
    \frac {dE_1} {dt}
    & =
    \left(
    \frac {1/2 -\varphi_1} {T_1}
    \frac {\partial T_1} {\partial E_1}
    +
    \frac{2\delta}{T_1}\int_0^{v_1}(\kappa^{-1})''(E_1-v^2)dv
    \right)
    (\varepsilon W \kappa_\delta' (X -x_1)).
\end {split}
\]
Using that $\kappa_\delta'(X-x_1) =\kappa' (\kappa^ {-1}
(E_1-v_1^2))/\delta=(\delta(\kappa^ {-1})' (E_1-v_1^2))^ {-1} $, we
find that
\[
    \frac {\partial \varphi_1} {\partial E_1}
    \frac {dE_1} {dt}
    =
    \varepsilon\mathcal{O}\left (\frac {1/2-\varphi_1}
    {\delta}\right)
    +
    \varepsilon\mathcal{O}
    \left(
    \frac {1}
    {(\kappa^ {-1})' (E_1-v_1^2)}
    \int_0^{v_1}(\kappa^{-1})''(E_1-v^2)dv
    \right).
\]
But here $1/2-\varphi_1$ is $\mathcal{O}(\delta)$.  See the proof of
Lemma \ref{lem:1d_smooth_delta} below.  Thus the claims about
$d\varphi_1/dt$ will be proven, provided we can uniformly bound
\[
    \frac{1}{(\kappa^ {-1})'(E_1-v_1^2)}
    \int_0^{v_1}(\kappa^{-1})''(E_1-v^2)dv.
\]
Note that the apparent divergence of the integral as
$\abs{v_1}\rightarrow\sqrt {E_1}  $ is entirely due to the fact that
our expression for $\varphi_1$ from Equation
\eqref{eq:anglevar_defn3} requires $\abs {v_1} < \sqrt E_1 $. If we
make the substitution $u=E_1-v^2$ and let $e=E_1-v_1^2$, then it
suffices to show that
\[
    \sup_{\mathcal{E}\leq E_1\leq\kappa(0)-\mathcal{E}}
    \;
    \sup_{0<e\leq E_1}
    \abs{\frac{1}{(\kappa^ {-1})' (e)}
    \int_e^{E_1}\frac{(\kappa^{-1})''(u)}{\sqrt{E_1-u}}du}
    <+\infty.
\]
The only difficulties occur when $e$ is close to $0$.  Thus it
suffices to show that
\[
    \sup_{\mathcal{E}\leq E_1\leq\kappa(0)-\mathcal{E}}
    \;
    \sup_{0<e\leq \mathcal{E}/2}
    \abs{\frac{1}{(\kappa^ {-1})' (e)}
    \int_e^{\mathcal{E}/2}\frac{(\kappa^{-1})''(u)}{\sqrt{E_1-u}}du}
\]
is finite.  But this is bounded by
\[
    \sup_{0<e\leq \mathcal{E}/2}
    \abs{\frac{1}{(\kappa^ {-1})' (e)}
    \int_e^{\mathcal{E}/2}\frac{(\kappa^{-1})''(u)}{\sqrt{\mathcal{E}/2}}du}
    =
    \sup_{0<e\leq \mathcal{E}/2}
    \abs{\frac{\sqrt{2/\mathcal{E}}}{(\kappa^ {-1})' (e)}
    \bigl((\kappa^ {-1})' (\mathcal{E}/2) -(\kappa^ {-1})' (e)\bigr) }
    ,
\]
which is finite because $(\kappa^ {-1})' (e)\rightarrow -\infty $ as
$e\rightarrow 0^ + $.  The claims about $d\varphi_2/dt$ can be
proven similarly.

\end {proof}

\paragraph* {Proof of Lemma \ref{lem:1d_smooth_delta}:}

\begin{proof}
We continue in the notation of the proofs of Lemmas
\ref{lem:1d_smooth_periods} and \ref{lem:1d_smooth_ode} above, and
we set $m_1 =2 $.  Then from Equation \eqref{eq:anglevar_defn3}, we
see that $\kappa_\delta' (X-x_1) =0 $ unless $\abs{\varphi_1
-1/2}\leq \abs{\frac{2\delta}{T_1} \int_0^{\sqrt
E_1}(\kappa^{-1})'(E_1-v^2)dv} = \delta F (E_1)/T_1 =\mathcal{O}
(\delta)  $.  Dealing with $\varphi_2 $ is similar.

\end {proof}

\textbf{Acknowledgments.}  The author is grateful to D. Dolgopyat,
who first introduced him to this problem, and who generously shared
his unpublished notes on averaging.  The author is also grateful to
 many people for encouragement and useful discussions regarding this
project, including G. Ariel, M. Lenci, K. Lin, and especially his
advisor, L.-S. Young. Thanks is due to an anonymous referee for
supplying several references relevant to the adiabatic piston
problem.  This research was partially supported by the National
Science Foundation Graduate Research Fellowship Program.

\bibliographystyle{amsplain}
\bibliography{pistonBib}
\end{document}